\documentclass[npg, manuscript]{copernicus}

\AtBeginDocument{\nolinenumbers}

\begin{document}

\title{A Continuous-Time Ensemble Kalman--Bucy Smoother for Causal Inference and Model Discovery}


\Author[1]{Zhang}{Jiang}
\Author[1]{Marios}{Andreou}
\Author[2]{Sebastian}{Reich}
\Author[1][chennan@math.wisc.edu]{Nan}{Chen}

\affil[1]{Department of Mathematics, University of Wisconsin-Madison, 480 Lincoln Drive, Madison, WI 53706, USA}
\affil[2]{Institut f\"ur Mathematik, University of Potsdam, Karl-Liebknecht-Str.~24/25, 14476 Potsdam, Germany}




\runningtitle{Continuous-time EnKBS for causal inference and model discovery}

\runningauthor{Z. Jiang et al.}

\received{}
\pubdiscuss{} 
\revised{}
\accepted{}
\published{}


\firstpage{1}

\maketitle

\begin{abstract}
    Data assimilation (DA) integrates observational information with model predictions to improve state estimation in complex systems. While filtering provides the basis for online forecasts by using only past and present observations, it can exhibit delays and biases when the underlying dynamics evolve rapidly or undergo regime transitions. Smoothing, which additionally incorporates future observations, provides a natural pipeline for hindcasting and reanalysis that yields an uncertainty reduction beyond the filter. This paper introduces an ensemble Kalman--Bucy smoother (EnKBS) for continuous-time DA of nonlinear dynamical systems, where the smoother's conditional distributions are reconstructed using ensemble moments. The result is a derivative-free framework that does not require explicit computation of tangent-linear or adjoint models, which recovers the exact smoothing mean and covariance equations in the infinite-ensemble limit for linear--Gaussian systems. Incorporating standard regularization techniques for high-dimensional systems, such as covariance localization and inflation, the skill of the EnKBS is demonstrated in various important scientific problems. By integrating future observations, which reveal the underlying causal mechanisms for retrospective state updates, the EnKBS is used for Bayesian-based inference of causal relationships and their temporal influence range in a dyadic trigger-feedback model and the development of a causality-driven iterative learning algorithm that identifies the structure and recovers the hidden parameters of a nonlinear reduced-order model mimicking midlatitude atmospheric circulation. Notably, both tasks remain effective with an ensemble size of $O(10)$ under partial observations, suggesting that EnKBS can support the instantaneous discovery of high-dimensional complex systems over time.
\end{abstract}


\introduction  
Complex dynamical systems arise in many scientific areas, such as atmospheric and ocean science, fluid dynamics, geophysics, and neuroscience \citep{chen2023stochastic,majda2016introduction,strogatz2024nonlinear,izhikevich2007dynamical,holmes2012turbulence}. Such systems are characterized by high dimensionality, strong nonlinearity, and nontrivial driving mechanisms \citep{strogatz2024nonlinear,kolmogorov1995turbulence,bar2019dynamics}.

A central mathematical problem when studying complex dynamical systems is reconstructing the full state of the system from partial and noisy observations, which is essential for system reanalysis, uncertainty quantification, forecasting, and optimal control \citep{chen2025taming,majda2012filtering,kalnay2003atmospheric,asch2016data,stengel1994optimal}. More importantly, accurate state estimation is often an essential prerequisite for model identification and the discovery of the underlying causal structure \citep{brunton2016discovering,chen2023causality,andreou2026assimilative}.

However, accurate state estimation is challenging. Uncertainties arise primarily from two sources: measurement noise due to imperfect observation tools or methodologies, and model error stemming from incomplete physical understanding or limited computational resolution. In chaotic systems, these uncertainties are amplified over time through multiscale nonlinear interactions, so state estimates that rely entirely on model simulations can quickly become unreliable \citep{kalnay2003atmospheric,lorenz1969predictability,lorenz2006predictability}. In contrast, purely data-driven approaches that rely on statistical interpolation of time-series data often fail to incorporate temporal dynamical knowledge and to enforce physics constraints \citep{chen2023causality,raissi2019physics}. A strategy is therefore needed that integrates observational information with model dynamics while addressing the limitations of each. Data assimilation (DA) has been the de facto approach, leading to a powerful framework that systematically integrates observations with dynamics to improve state estimation and reduce uncertainty and biases \citep{law2015data,evensen2009data,asch2016data}. Over the past few decades, DA and its variants have been extensively developed and deployed across various disciplines, with notable successes in numerical weather prediction and medical imaging \citep{bauer2015quiet,lorenc1986analysis,tarantola2005inverse}.

Broadly speaking, classical DA methods are categorized into filters and smoothers. Filters estimate the state using observations up to and including the present, whereas smoothers incorporate information from the entire observation window, including past, present, and future observations. While filtering is the standard approach for real-time forecasting, smoothing is theoretically more accurate for retrospective analysis \citep{simon2006optimal,evensen2000ensemble,buehner2017ensemble}. On the one hand, smoothing systematically improves estimation error and reduces uncertainty \citep{simon2006optimal,moore1979optimal}. On the other hand, smoothing can mitigate the detection delay that filters often exhibit at the onset of extreme events \citep{andreou2026assimilative,andreou2024adaptive,moser2026mechanisms}
(see also Sect.~\ref{subsec:ACI}). This distinction becomes particularly important in systems where the underlying dynamics are driven by hidden forcing mechanisms. As a result, filters tend to detect the driving signal only after its effects have been observed \citep{chen2016filtering}. In contrast, by leveraging future information, smoothers can more accurately identify the onset and persistence of the underlying driving mechanisms \citep{andreou2026assimilative,andreou2024adaptive,andreou2025bridging}.

Although a wide range of DA methods has been developed, many methods face theoretical and practical limitations. The Kalman filter and the Rauch--Tung--Striebel (RTS) smoother provide optimal solutions, but only for linear--Gaussian systems \citep{kalman1960new,rauch1965maximum}. In complex dynamical systems, however, strongly non-Gaussian features such as skewness and heavy tails often violate the Gaussian assumptions underlying the classical theory \citep{jazwinski2007stochastic}. For general nonlinear dynamics, variational DA approaches typically provide deterministic state estimates that do not quantify uncertainty, while particle methods suffer from severe computational burdens in high dimensions due to the curse of dimensionality \citep{snyder2008obstacles,rabier2005overview,lorenc2003potential}. Consequently, ensemble Kalman methods have emerged as an effective alternative \citep{evensen1994sequential,burgers1998analysis}. By approximating system statistics with a finite ensemble, these methods provide state estimates that explicitly quantify uncertainty while remaining derivative-free and computationally feasible in many high-dimensional settings. In practice, the ensemble size (typically $O(10^2)$) is often chosen independently of the system's effective dimension \citep{evensen2003ensemble,kalnay20074}. Due to these advantages, ensemble methods have been widely adopted in large-scale operational applications, including numerical weather prediction, ocean science, geophysics, and fluid dynamics \citep{houtekamer2016review,iglesias2013ensemble,evensen2003ensemble,aanonsen2009ensemble}. Most existing ensemble filtering and smoothing methods, however, are formulated in discrete time, reflecting the sequential structure of forecast-analysis cycles. Efficient continuous-time smoothers for nonlinear systems remain comparatively less explored, despite their importance for systems naturally described by stochastic
differential equations. Existing formulations, such as the continuous-time unscented RTS smoother \citep{sarkka2010continuous}, directly propagate the smoother mean and covariance. In contrast, our ensemble-based smoother evolves stochastic trajectories, providing conditional samples for downstream tasks such as model discovery, while naturally admitting covariance localization for high-dimensional systems.

In this paper, we develop a continuous-time ensemble Kalman--Bucy smoother (EnKBS) for efficient state estimation and trajectory sampling in complex nonlinear dynamical systems. To address high-dimensional settings, the method incorporates standard regularization techniques, including covariance localization and inflation, which mitigate sampling errors, improve numerical stability, and prevent ensemble divergence in nonlinear models \citep{anderson2001ensemble,bergemann2010localization,houtekamer1998data}. Beyond state estimation, we demonstrate that EnKBS provides a powerful tool for dynamical discovery. We incorporate EnKBS into a Bayesian causal inference framework that identifies causal relationships by tracing driving signals backwards from their observed effects \citep{andreou2026assimilative}. Numerical tests are presented using a nonlinear dyad model exhibiting extreme events, for which closed-form filtering and smoothing distributions are available \citep{andreou2026assimilative,andreou2024martingale}. We verify asymptotic convergence of the method to the exact analytical solution and show that it remains accurate and statistically consistent even with an ensemble size of $O(10)$. We further integrate EnKBS into a causality-based iterative learning algorithm for partially observed systems, where EnKBS-based trajectory sampling enables the recovery of hidden parameters and identification of the underlying model structure \citep{chen2023causality}.

The rest of this paper is organized as follows. Sect.~\ref{sec:EnKBS} formulates the problem and introduces the proposed continuous-time EnKBS, together with regularization strategies for high-dimensional systems. Sect.~\ref{sec:cases} evaluates EnKBS on the high-dimensional Lorenz-96 system \citep{lorenz2006predictability} and presents applications to assimilative causal inference (ACI) and causality-based model discovery \citep{chen2023causality, andreou2026assimilative}. We conclude the paper with a discussion in Sect.~\ref{sec:Con}. Mathematical derivations of the EnKBS and its connection to the discrete-time RTS smoother in a formal time limit are deferred to Appendix~\ref{app:der}. Finally, Appendix~\ref{app:moment_consistency} shows that the proposed EnKBS recovers the exact mean and covariance equations of the classical continuous-time RTS smoother in the linear--Gaussian setting.

\section{The Continuous-Time Ensemble Kalman--Bucy Smoother}\label{sec:EnKBS}

\subsection{Problem Setup}\label{subsec:problem}
 We consider a general class of continuous-time nonlinear systems. Let $\vec{x}(t) \in \mathbb{R}^{n_{\vec{x}}}$ denote the hidden states and $\vec{y}(t)\in \mathbb{R}^{n_{\vec{y}}}$ denote the observations (data). The dynamics are governed by the following stochastic differential equations (SDEs) \citep{jazwinski2007stochastic,sarkka2019applied}:
 \begin{subequations}\label{eq:general_system}
\begin{align}
 \mathrm{d}\vec{x} &=\vec{f}(\vec{x}, \vec{y}, t)\mathrm{d} t + \mathbf \Sigma(\vec{y}, t)^{1/2} \mathrm{d}\vec{B}(t), \label{eq:signal_process} \\
\mathrm{d}\vec{y} &= \vec{h}(\vec{x}, \vec{y}, t)\mathrm{d} t + \mathbf \Gamma(\vec{y}, t)^{1/2} \mathrm{d}\vec{W}(t). \label{eq:obs_process}
\end{align}
\end{subequations}
$\vec{f}$ is the drift of the signal process Eq.~(\ref{eq:signal_process}) and $\vec{h}$ is the drift of the observation process Eq.~(\ref{eq:obs_process}). The terms $\vec{B}(t)\in \mathbb{R}^{n_{\vec{x}}}$ and $\vec{W}(t) \in \mathbb{R}^{n_{\vec{y}}}$ are independent standard Wiener processes. For simplicity, we assume that the signal process noise covariance $\mathbf\Sigma(\vec{y},t)\in\mathbb{R}^{n_{\vec{x}}\times n_{\vec{x}}}$ and the observation noise covariance $\mathbf\Gamma(\vec{y},t)\in\mathbb{R}^{n_{\vec{y}}\times n_{\vec{y}}}$ are symmetric and positive definite. Throughout the paper, all SDEs, whether written in differential form or in formal derivative form, are understood in the It\^o sense unless otherwise stated \citep{jazwinski2007stochastic,oksendal2003Stochastic}.

We consider a fixed time horizon $t\in[0,T]$, and aim to characterize the conditional distribution of $\vec{x}(t)$ given the available observations. Let $\vec{y}_{[0, \tau]} = \{ \vec{y}(s) : 0 \le s \le \tau \}$ denote the observation trajectory up to time $\tau\in[0,T]$. We formulate the filtering and smoothing problems of DA as follows:
\begin{enumerate}
    \item \textbf{Filtering:} Estimate the conditional distribution of $\vec{x}(t)$ given the observational history up to the present $\vec{y}_{[0,t]}$. The filtering distribution is denoted by $p(\vec{x}(t)\mid \vec{y}_{[0,t]})$, with mean $\overline{\vec{x}}_{\mathrm{f}}(t)$ and covariance $\mathbf P_{\mathrm{f}}(t)$.
    \item \textbf{Smoothing:} Estimate the conditional distribution of $\vec{x}(t)$ given the full observation trajectory $\vec{y}_{[0,T]}$, including the past and the future. The smoothing distribution is denoted by $p(\vec{x}(t)\mid \vec{y}_{[0,T]})$, with mean $\overline{\vec{x}}_{\mathrm{s}}(t)$ and covariance $\mathbf P_{\mathrm{s}}(t)$.
\end{enumerate}

It is worth mentioning that the coefficients in~Eq.~(\ref{eq:general_system}), including the drifts $\vec{f}$ and $\vec{h}$, as well as the noise covariances $\mathbf\Sigma$ and $\mathbf\Gamma$, may depend explicitly on the observation history $\vec{y}_{[0,t]}$. Since $\vec{y}$ is given in the filter and smoother, its trajectory can be treated as a time-varying parameter \citep{crisan2011oxford,rozovsky2018stochastic}. This flexibility allows for broad applications in which the observed variables actively feed back on and modulate the hidden dynamics; however, for our purposes, we assume the Markovian structure in Eq.~(\ref{eq:general_system}) \citep{oksendal2003Stochastic}. For notational simplicity, we may omit the explicit dependence of the model components on the dependent and time variables when the context is clear, and write $\vec{f}(\vec{x},\vec{y},t)$ as $\vec{f}$, and similarly for $\vec{h}$, $\mathbf\Sigma$, and $\mathbf\Gamma$.

In general, evolving the full filtering or smoothing distribution is computationally infeasible \citep{crisan2011oxford}. Therefore, we focus on estimating the first two moments: the mean and covariance. The idea is that in systems with linear dynamics and Gaussian noise, the conditional distributions are Gaussian and fully described by their mean and covariance, i.e., \(p(\vec{x}(t) \mid \vec{y}_{[0, t]})=\mathcal N(\overline{\vec{x}}_{\mathrm{f}}(t),\mathbf P_{\mathrm{f}}(t))\) and  \(p(\vec{x}(t) \mid \vec{y}_{[0, T]})=\mathcal N(\overline{\vec{x}}_{\mathrm{s}}(t),\mathbf P_{\mathrm{s}}(t))\). These exact moments are governed by the classical Kalman--Bucy Filter \citep{kalman1961new} and continuous-time Rauch--Tung--Striebel (RTS) smoother equations \citep{rauch1965maximum}.

For general nonlinear systems where closed-form DA solutions are unavailable, we adopt an ensemble-based approach \citep{evensen2009data,evensen1994sequential}. Here, the conditional distributions are represented by a finite ensemble of size $m$. This can be viewed as a Monte Carlo approximation in which the conditional moments are replaced by empirical ensemble moments. Each ensemble member independently evolves according to Eq.~(\ref{eq:general_system}), but the ensemble constituents are coupled through the empirical mean and covariance. Since individual ensemble anaylsis members are propagated through the nonlinear
dynamics, the method is not restricted to linear--Gaussian systems, and
the empirical ensemble distribution can retain certain non-Gaussian features.
However, because ensemble methods rely on empirical
first- and second-order statistics, accurate representation of general non-Gaussian distributions is not guaranteed.

Let $\vec{x}^{(i)}(t)$, $i=1,\dots,m$, denote the ensemble members. The empirical mean $\overline{\vec{x}}(t)$ and covariance $\mathbf P(t)$ for $t\in[0,T]$ are defined by
    \begin{equation*}
        \overline{\vec{x}}(t):=\frac{1}{m}\sum_{i=1}^m \vec{x}^{(i)}(t),\qquad
        \mathbf P(t):=\frac{1}{m-1}\sum_{i=1}^m\big(\vec{x}^{(i)}(t)-\overline{\vec{x}}(t)\big)\big(\vec{x}^{(i)}(t)-\overline{\vec{x}}(t)\big)^\top.
    \end{equation*}
When needed, we use subscripts $\mathrm{f}$ and $\mathrm{s}$ to distinguish between the filter- and smoother-based ensembles. For brevity, we use $\big(\overline{\vec{x}}_{\mathrm{f}},\mathbf P_{\mathrm{f}}\big)$ and $\big(\overline{\vec{x}}_{\mathrm{s}},\mathbf P_{\mathrm{s}}\big)$ to denote the corresponding empirical moments, and do not distinguish them from the true moments unless confusion may arise.

While filtering is the standard approach for online state estimation,
when the full observation trajectory over $[0,T]$,
$\vec{y}_{[0,T]}$, is available, applications such as Earth system
reanalysis \citep{kalnay2003atmospheric} and offline model and parameter
learning \citep{navon1998practical,yang2009using} can benefit from
retrospective analysis. Based on the filtering quantities $\big(\overline{\vec{x}}_{\mathrm{f}}(t), \mathbf P_{\mathrm{f}}(t)\big)$, the smoothing quantities $\big(\overline{\vec{x}}_{\mathrm{s}}(t), \mathbf P_{\mathrm{s}}(t)\big)$ incorporate information from future observations $\vec{y}_{[t,T]}$, producing more accurate estimates and reducing uncertainty over $[0,T]$. This improvement motivates a systematic forward--backward architecture \citep{evensen2009data,sarkka2023bayesian}:
    \begin{enumerate}
        \item \textbf{Forward pass (filtering):} Integrate forward from $t=0$ to $t=T$ to compute the filter ensemble \(\vec{x}^{(i)}_{\mathrm{f}}(t),i=1,\dots,m\) (forecasting). The filtering mean and covariance $\{ \overline{\vec{x}}_{\mathrm{f}}(t), \mathbf P_{\mathrm{f}}(t)\}_{0\le t\le T}$ are obtained.
        \item \textbf{Backward pass (smoothing):} Integrate backward-in-time from $t=T$ to $t=0$. Using the filtering estimates as a prior, this pass produces the smoothing ensemble \(\vec{x}_{\mathrm{s}}^{(i)}(t),i=1,\dots,m\) and therefore its statistics $\{ \overline{\vec{x}}_{\mathrm{s}}(t),\mathbf P_{\mathrm{s}}(t)\}_{0\le t\le T}$ (hindcasting).
    \end{enumerate}

In our method, we implement this architecture by utilizing the continuous-time ensemble Kalman--Bucy filter (EnKBF) for the forward pass \citep{bergemann2012Ensemble}, and the EnKBS for the backward pass. The operational details of these methods are presented in the following subsections.

\subsection{The Forward Pass: EnKBF} \label{subsec:EnKBF}
We present a stochastic variant of the EnKBF in this subsection, dropping the subscripts \(\mathrm{f}\) for notational simplicity. We first consider a linear observation drift
\(\vec{h}(\vec{x},t)=\mathbf H\vec{x}\).
For each ensemble member, we generate a simulated observation increment
\begin{equation}\label{eq:sim-data}
    \mathrm{d} \vec{y}^{(i)}(t)
    := \mathbf H\vec{x}^{(i)}(t)\mathrm{d} t
    + \mathbf\Gamma^{1/2}\mathrm{d} \vec{W}^{(i)}(t),
\end{equation}
where $\{\vec{W}^{(i)}\}_{i=1,\ldots,m}$ are mutually independent standard Wiener processes.
For a linear observation operator \(\mathbf H\), the stochastic EnKBF is given by \citep{bergemann2012Ensemble}:
\begin{equation}\label{eq:enkbf-sde-linear}
    \mathrm{d} \vec{x}^{(i)}
    = \vec{f}(\vec{x}^{(i)},t)\mathrm{d} t
    + \mathbf \Sigma^{1/2}\mathrm{d}\vec{B}^{(i)}
    + \mathbf P \mathbf H^\top \mathbf\Gamma^{-1}
    \left(\mathrm{d} \vec{y} - \mathrm{d} \vec{y}^{(i)}\right).
\end{equation}
A formal connection between this continuous-time formulation and the standard discrete-time ensemble Kalman filter is provided in Appendix~\ref{app:der}. 
Treating the simulated observations as random variables corresponding to each ensemble member ensures statistical consistency with the Riccati covariance equation in classical Kalman filter theory; in the absence of the noise term $\mathbf\Gamma^{1/2}\mathrm{d} \vec{W}^{(i)}$ in Eq.~(\ref{eq:sim-data}), the filtering covariance would be systematically underestimated
\citep{burgers1998analysis,sakov2008deterministic}. 

In practical applications, the observation drift $\vec{h}$ is often nonlinear. To extend Eq.~(\ref{eq:enkbf-sde-linear}) to this setting without explicit Jacobian computations, we use the following ensemble approximations:
\begin{gather}
    \mathbf H\vec{x}^{(i)} \approx \vec{h}(\vec{x}^{(i)}, t),\label{eq:abb_tl}\\
    \mathbf P\mathbf H^\top \approx \mathbf P_{\vec{x}\vec{h}}.\label{eq:abb_aj}
\end{gather}
Here, $\mathbf P_{\vec{x}\vec{h}}$ denotes the cross-covariance matrix between the ensemble and its observation-operator-evolved counterpart:
\begin{equation*}
     \mathbf P_{\vec{x}\vec{h}}:=\frac{1}{m-1}\sum_{i=1}^m\left(\vec{x}^{(i)}-\overline {\vec{x}}\right)\left(\vec{h}(\vec{x}^{(i)})-\overline{\vec{h}}\right)^\top,
\end{equation*}
where
\[
\overline{ \vec{h}}(t):=\frac{1}{m}\sum_{i=1}^m \vec{h}(\vec{x}^{(i)},t).
\]
Substituting Eq.~(\ref{eq:abb_tl})--Eq.~(\ref{eq:abb_aj}) into Eq.~(\ref{eq:enkbf-sde-linear}) yields the practical derivative-free EnKBF formulation for \(i=1,\dots,m\):
\begin{equation}\label{eq:enkbf-sde-stc}
    \mathrm{d} \vec{x}^{(i)}
    = \vec{f}(\vec{x}^{(i)},t)\mathrm{d} t
    + \mathbf \Sigma^{1/2}\mathrm{d}\vec{B}^{(i)}
    +\mathbf P_{\vec{x}\vec{h}}\mathbf\Gamma^{-1}\left(\mathrm{d} \vec{y} - \mathrm{d} \vec{y}^{(i)}\right),
\end{equation}
where, under the approximation (\ref{eq:abb_tl}), the simulated observation in Eq.~(\ref{eq:sim-data}) becomes $\mathrm{d} \vec{y}^{(i)} = \vec{h}(\vec{x}^{(i)},t)\mathrm{d} t + \mathbf\Gamma^{1/2}\mathrm{d} \vec{W}^{(i)}$.
We remark that applying Eq.~(\ref{eq:abb_tl})--Eq.~(\ref{eq:abb_aj}) avoids explicit tangent-linear or adjoint computations related to \(\mathbf H\) \citep{asch2016data,evensen2003ensemble}. The key assumption behind Eq.~(\ref{eq:abb_aj}) is
\[
\vec{h}({\vec{x}}^{(i)}) - \overline{\vec{h}}
\approx
\mathbf H\big({\vec{x}}^{(i)} - \overline {\vec{x}}\big),
\]
which can be interpreted as a finite-difference approximation of the drift $\vec{h}$ when the $i$-th ensemble spread, ${\vec{x}}^{(i)} - \overline {\vec{x}}$, is sufficiently small.

\subsection{The Backward Pass: EnKBS} \label{subsec:EnKBS}

Given a filtering ensemble \( \vec{x}_{\mathrm{f}}^{(i)}(t),i=1,\dots,m,\) produced by the forward stochastic EnKBF Eq.~(\ref{eq:enkbf-sde-stc}), the smoother integrates backward from \(t=T\) to \(t=0\) to refine the ensemble statistics over $[0,T]$ in the context of the full observed trajectory $\vec{y}_{[0,T]}$. Denote the smoothing ensemble by \( \vec{x}_{\mathrm{s}}^{(i)}(t),i=1,\dots,m\). Initializing the backward pass at $t=T$ by
\begin{equation}
    \vec{x}_{\mathrm{s}}^{(i)}(T) = \vec{x}_{\mathrm{f}}^{(i)}(T),
\end{equation}
the backward-in-time EnKBS ensemble member dynamics for \( i=1,\dots,m\) are given by
\begin{equation}\label{eq:enkbs_member_sde}
    \mathrm{d}\vec{x}_{\mathrm{s}}^{(i)}
    =
    \vec{f}\big(\vec{x}_{\mathrm{s}}^{(i)},t\big)\mathrm{d}t
    +
    \mathbf\Sigma^{1/2}\mathrm{d}\vec{B}^{(i)}
    +
    \mathbf\Sigma \mathbf P_{\mathrm{f}}^{-1}
    \big(\vec{x}_{\mathrm{s}}^{(i)}-\vec{x}_{\mathrm{f}}^{(i)}\big)
    \mathrm{d}t.
\end{equation}
Here, \(\mathbf P_{\mathrm{f}}\) is the empirical covariance of the filtering ensemble \( \vec{x}_{\mathrm{f}}^{(i)}\). In Eq.~(\ref{eq:enkbs_member_sde}), we use the same Wiener process paths $\vec{B}^{(i)}$ as in the forward filter Eq.~(\ref{eq:enkbf-sde-stc}), as justified by the derivations in Appendix~\ref{app:der}.

We next specify the Euler--Maruyama time discretization used to evaluate
Eq.~(\ref{eq:enkbs_member_sde}) backward from \(t=T\) to \(t=0\). Let $\tau>0$ denote the time step and $t_k = k\tau$. We use the subscript $(\cdot)_k$ to denote evaluation at time $t_k$, e.g., $\vec{x}_k=\vec{x}(t_k)$. The forward EnKBF Eq.~(\ref{eq:enkbf-sde-stc}) is discretized as
\begin{equation}\label{eq:enkbf-disc}
    \vec{x}_{\mathrm{f},k+1}^{(i)}
    = \vec{x}_{\mathrm{f},k}^{(i)}
    + \tau \vec{f}\left(\vec{x}_{\mathrm{f},k}^{(i)},t_k\right)
    + \sqrt{\tau}\mathbf\Sigma_k^{1/2}\vec{B}_k^{(i)}
    + \mathbf P_{\vec{x}\vec{h},k}\mathbf\Gamma^{-1}\left(\Delta\vec{y}_{k+1} - \Delta \vec{y}_{k+1}^{(i)}\right),
\end{equation}
for $k=0,1,\dots,K-1$, where $\vec{B}_k^{(i)},\vec{W}_k^{(i)}
\sim\mathcal{N}(\mathbf 0,\mathbf{I})$ are independent, $\Delta\vec{y}_{k+1}=\vec{y}_{k+1}-\vec{y}_k$ is the observation increment, and
\begin{equation*}
    \Delta \vec{y}_{k+1}^{(i)} := \tau\vec{h}\left(\vec{x}_{\mathrm{f},k}^{(i)},t_k\right) + \sqrt{\tau}\mathbf\Gamma^{1/2}\vec{W}_k^{(i)}
\end{equation*}
is the simulated observation increment for the $i$-th member, consistent with the continuous-time definition Eq.~(\ref{eq:sim-data}). The backward EnKBS Eq.~(\ref{eq:enkbs_member_sde}) is discretized as
\begin{equation}\label{eq:enkbs-disc}
    \vec{x}_{\mathrm{s},k}^{(i)}
    = \vec{x}_{\mathrm{s},k+1}^{(i)}
    - \tau \vec{f}\left(\vec{x}_{\mathrm{s},k+1}^{(i)},t_{k+1}\right)
    - \sqrt{\tau}\mathbf\Sigma_k^{1/2}\vec{B}_k^{(i)}
    - \tau\mathbf\Sigma_{k+1}\mathbf P_{\mathrm{f},k+1}^{-1}
      \left(\vec{x}_{\mathrm{s},k+1}^{(i)}
      - \vec{x}_{\mathrm{f},k+1}^{(i)}\right),
\end{equation}
for $k=K-1,K-2,\dots,0$, with initialization $\vec{x}_{\mathrm{s},K}^{(i)}
= \vec{x}_{\mathrm{f},K}^{(i)}$. The backward EnKBS is advanced explicitly from \(t_{k+1}\) to \(t_k\),
with the deterministic terms evaluated at \(t_{k+1}\). The stochastic
increment over \([t_k,t_{k+1}]\) is retained from the forward pass, so
the same \(\vec{B}_k^{(i)}\), rather than \(\vec{B}_{k+1}^{(i)}\), is
used in the backward step.

As a remark, in the EnKBS formulation Eq.~(\ref{eq:enkbs_member_sde}), the backward dynamics do not depend explicitly on the observation increment $\mathrm{d}\vec{y}$ as with the EnKBF in Eq.~(\ref{eq:enkbf-sde-stc}). The implicit dependence on $\vec{y}$ is instead via the filtering ensemble ${\vec{x}_{\mathrm{f}}^{(i)}(t)}$ computed in the forward pass. Under the assumption of independent signal and observation white noises, the information in $\vec{y}$ is already encoded in the future state $\vec{x}$ and propagates backward to correct the current estimate along the smoother pass. Conditioned on the future \(\vec{x}\), the Markov property of the system yields that including any observation \(\vec{y}\) depending on it does not provide additional information in the absence of cross-interacting noise \citep{andreou2024martingale, crisan2011oxford, durrett2019probability}.

Detailed derivations of Eq.~(\ref{eq:enkbs_member_sde}) and its formal connection to the discrete-time RTS smoother in the formal limit of a vanishing observation interval are provided in Appendix~\ref{app:der}. In addition, Appendix~\ref{app:moment_consistency} shows that the proposed EnKBS recovers the exact mean and covariance equations of the continuous-time RTS smoother in the linear--Gaussian setting \citep{rauch1965maximum}.

\subsection{Localization and Inflation} \label{subsec:LocInf}

In high-dimensional data assimilation with state dimension $n$, a small ensemble size $m\ll n$ inevitably incurs sampling error and cannot fully resolve all dynamic modes. Specifically, the full covariance matrix is $n\times n$, whereas the empirical ensemble covariance $\mathbf P$ is rank-deficient with rank at most $m-1$, which is often insufficient to fully capture the dominant unstable modes of the system \citep{bocquet2017four}. This rank deficiency may cause instability or divergence of the filter and smoother and induce spurious, nonphysical long-range correlations \citep{evensen2003ensemble,houtekamer1998data}. In many operational applications, affordable ensemble sizes typically range from $O(10)$ to $O(100)$, which can be much smaller than the dimension of the unstable subspace \citep{asch2016data}. Increasing $m$ to that scale is often computationally prohibitive.

Two standard regularization techniques, localization and inflation, are widely used to address these issues and improve stability in high-dimensional systems. These ad hoc strategies, to some extent, help explain the practical success of ensemble methods in applications such as numerical weather prediction and the geosciences \citep{houtekamer2016review}.

\noindent\textbf{Localization.} Covariance localization, also called Schur localization, regularizes the empirical covariance by damping long-range state-space correlations and thereby increases its rank \citep{asch2016data, bergemann2010localization}. A standard approach applies a Schur (element-wise) product between the empirical covariance $\mathbf P$ and a short-ranged correlation matrix $\mathbf C_{\mathrm{loc}}\in\mathbb{R}^{n\times n}$:
    \begin{equation}\label{eq:Ploc-schur}
        \mathbf P\leftarrow\mathbf C_{\mathrm{loc}}\circ \mathbf P,
    \end{equation}
where $\circ$ denotes the Schur product of matrices.

By the Schur product theorem \citep{horn1985matrix}, $\mathbf C_{\mathrm{loc}}\circ \mathbf P$ is positive (semi-)definite whenever $\mathbf C_{\mathrm{loc}}$ and $\mathbf P$ are positive (semi-)definite. A common choice constructs $\mathbf C_{\mathrm{loc}}$ from the Gaspari--Cohn fifth-order piecewise rational function $G$ \citep{gaspari1999construction},
\begin{equation*}
    G(r)=
    \begin{cases}
        1 - \dfrac{5}{3}r^2 + \dfrac{5}{8}r^3 + \dfrac{1}{2}r^4 - \dfrac{1}{4}r^5,
        & 0 \le r < 1,\\
        4 - 5r + \dfrac{5}{3}r^2 + \dfrac{5}{8}r^3 - \dfrac{1}{2}r^4 + \dfrac{1}{12}r^5 - \dfrac{2}{3r},
        & 1 \le r< 2,\\
        0, & r \ge 2,
    \end{cases}
\end{equation*}
for $r\ge 0$. It resembles a Gaussian distribution, while remaining compactly supported and inexpensive to evaluate. Given a localization radius $r_0>0$ and a distance $d(i,j)$ between
the indices $i$ and $j$ of the state components, we define
    \begin{equation*}
        (\mathbf C_{\mathrm{loc}})_{i,j} := G\left(\frac{d(i,j)}{r_0}\right),
    \end{equation*}
so that correlations are truncated beyond $d(i,j)>2r_0$.

Accordingly, the localized formulation of the stochastic EnKBF Eq.~(\ref{eq:enkbf-sde-stc}) replaces covariance \(\mathbf P_{\vec{x}\vec{h}}\) by \(\mathbf C_{\mathrm{loc,1}}\circ \mathbf P_{\vec{x}\vec{h}}\) \citep{khare2008investigation}: 
\begin{equation}\label{eq:enkbf-sde-stc_loc}
    \mathrm{d} \vec{x}^{(i)}
    = \vec{f}(\vec{x}^{(i)},t)\mathrm{d} t
    + \mathbf \Sigma^{1/2}\mathrm{d}\vec{B}^{(i)}
    +(\mathbf C_{\mathrm{loc,1}}\circ\mathbf P_{\vec{x}\vec{h}})\mathbf\Gamma^{-1}\left(\mathrm{d} \vec{y} - \mathrm{d} \vec{y}^{(i)}\right),
\end{equation}
and, for consistency, we also apply a localization to the backward EnKBS Eq.~(\ref{eq:enkbs_member_sde}):
\begin{equation}\label{eq:enkbs_member_sde_loc}
   \mathrm{d}\vec{x}_{\mathrm{s}}^{(i)}
=
\vec f(\vec{x}_{\mathrm{s}}^{(i)},t)\mathrm dt
+
\mathbf\Sigma^{1/2}\mathrm d\vec B^{(i)}
+
\mathbf\Sigma
(\mathbf C_{\mathrm{loc,2}}\circ\mathbf P_{\mathrm f})^{-1}
(\vec{x}_{\mathrm s}^{(i)}-\vec{x}_{\mathrm f}^{(i)})
\mathrm dt.
\end{equation}

\noindent\textbf{Inflation.} With localization, sampling error and spurious correlations are reduced, but still may accumulate over time. A practical solution to this is multiplicative inflation, which rescales ensemble anomalies around the mean by a factor $\delta^2>1$  \citep{anderson1999monte}:
    \begin{equation}
        \vec{x}^{(i)} \longleftarrow \overline{\vec{x}} + \delta\big(\vec{x}^{(i)}-\overline{\vec{x}}\big), \qquad i=1,\dots,m.
    \label{eq:inflation-members}
\end{equation}

The multiplicative inflation Eq.~(\ref{eq:inflation-members}) aims to compensate for sampling errors that are the indirect consequence of the nonlinearity in the ensemble forecast, which helps to cure the intrinsic error of ensemble filtering methods \citep{asch2016data}. In addition, inflation should be applied only for the forward filtering, not the backward smoothing pass. Since sampling errors are already accounted for in the filter, applying inflation in the backward update would result in a suboptimal solution \citep{cosme2010implementation}.

\section{EnKBS for State Estimation, Causal Inference, and Model Identification}\label{sec:cases}
In this section, the EnKBS is adopted to study three topics in complex nonlinear dynamical systems. The first study is the state estimation of a relatively high-dimensional chaotic system, which highlights the necessity of regularization and demonstrates the applicability of EnKBS in such systems. The remaining two examples are nonlinear coarse-grained conceptual models, where EnKBS is adopted for advancing causal inference and discovering the hidden model structure and parameters.

In all experiments, we use the Euler--Maruyama scheme for numerical time integration \citep{kloeden1992numerical}. For the backward smoothing pass, we discretize the formulation Eq.~(\ref{eq:enkbs_member_sde}) and propagate the ensemble backward in time. For all case studies, the reference trajectory is generated by
integrating the full stochastic dynamical system. The trajectories of
the designated observed variables, evaluated on the numerical time
grid, are used directly as the observations, without introducing an
additional discrete measurement model. The specific hidden--observed
variable partition and numerical time step are given in the
corresponding subsections.

\subsection{State Estimation in a Relatively High-Dimensional Chaotic System: The Lorenz-96 Model}\label{subsec:L96}

The first example considers the state estimation of the Lorenz-96 system. Introduced by Edward Lorenz, it was originally designed to mimic the propagation of atmospheric waves on a midlatitude zonal circle \citep{lorenz2006predictability,lorenz1998optimal}. It is now widely used as a benchmark test case in data assimilation. Let $\vec{x}=(x_1,\dots,x_n)^\top\in\mathbb{R}^n$.
We consider the stochastic Lorenz-96 dynamics:
    \begin{equation}\label{eq:L96}
        \dot  x_j = (x_{j+1}-x_{j-2})x_{j-1}-x_j+F+\sigma_j \dot  W_j,
        \qquad j=1,\dots,n.
    \end{equation}
Periodic indexing modulo $n$ is applied, i.e., $x_{j+n}=x_j$ for all integers $j$. Here, $W_j$ are mutually independent Wiener processes and $\dot{x}$ denotes derivative with respect to time \(t\).

The quadratic term \((x_{j+1}-x_{j-2})x_{j-1}\) mimics nonlinear advection and preserves the total energy, which is proportional to \(\sum_ix_i^2\). The linear term $-x_j$ represents mechanical or thermal dissipation and provides system damping, and the constant forcing $F$ drives the system into a chaotic regime when it exceeds a certain threshold value. With this parameter choice, the deterministic system has $13$ positive and $1$ neutral Lyapunov exponent, giving a $14$-dimensional unstable--neutral subspace. \citep{bocquet2017four, lorenz1998optimal}.

\subsubsection{The Experiment Setup}
In this test, we observe every other state component (the even indices) and treat the remaining components as hidden \citep{yang2009using}. In other words, the hidden signal process Eq.~(\ref{eq:signal_process}) and observational process Eq.~(\ref{eq:obs_process}) correspond respectively to
\begin{align*}
    \dot x_j &= (x_{j+1}-x_{j-2})x_{j-1}-x_j+F+\sigma_{\mathrm h}\dot W_j,
        \qquad j=1,3,5,\dots,n-1,\\
    \dot x_j& = (x_{j+1}-x_{j-2})x_{j-1}-x_j+F+\sigma_{\mathrm o} \dot  W_j,
        \qquad j=2,4,6,\dots,n.
\end{align*}
We set the noise variance in the observed processes to \(\sigma_{\mathrm o}^2=0.1\), and the noise variance in the unobserved processes as \(\sigma_{\mathrm h}^2=5\). The dimension of the system is $n=40$ and the forcing $F=8$ represents a strongly chaotic regime.

\begin{figure}[t]
    \includegraphics[width=8.3cm]{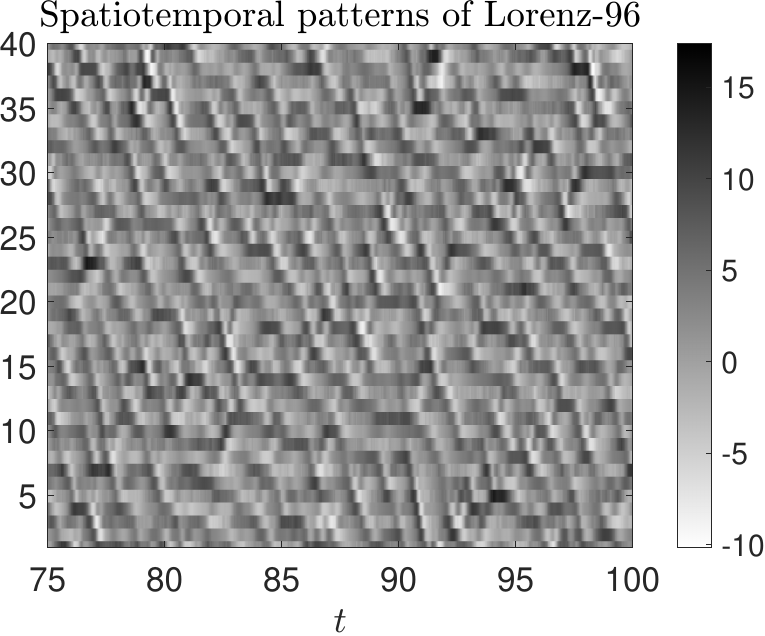}
    \caption{Lorenz-96 reference trajectory over $t\in[75,100]$. The horizontal axis is time, the vertical axis is the state index $j\in\{1,\dots,40\}$, and the color indicates the value of $x_j(t)$.}
    \label{fig:ref_l96}
\end{figure}

We set the observation rate to be equal to the numerical integration time step $\Delta t=0.0005$. A reference trajectory $\vec{x}_{\mathrm{ref}}$ is generated over the time window $[0,T]$ with $T=100$, and a segment is shown in Figure~\ref{fig:ref_l96}. We initialize the ensemble by adding small perturbations to the initial reference state. For localization, we use the periodic distance between indices $i$ and $j$:
    \begin{equation*}
        d(i,j):=\min\big\{|i-j|, n-|i-j|\big\},
        \qquad 1\le i,j\le n.
    \end{equation*}

\subsubsection{Results}
\begin{figure*}[t]
    \includegraphics[width=12cm]{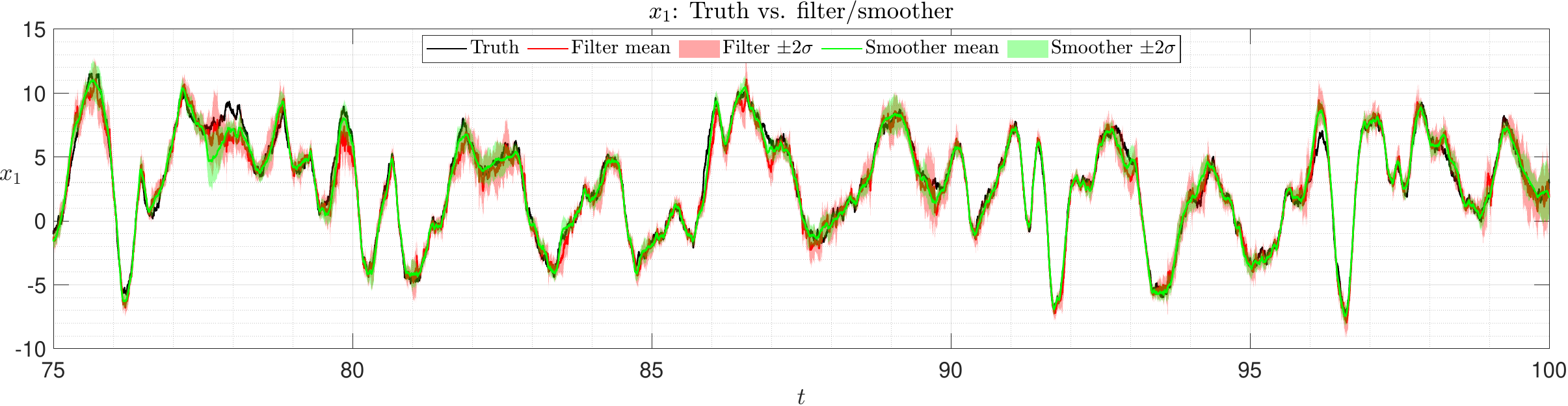}
    \caption{Trajectories of $x_1(t)$ over $t\in[75,100]$ for Lorenz-96. The black curve is the reference solution; the red and green curves are the EnKBF and EnKBS ensemble means, respectively. Shaded bands show the respective $\pm2$ standard-deviation intervals. Parameters: $m=10$, $r_0=4$, and $\delta^2=1.01$.}
    \label{fig:trj_l96}
\end{figure*}

The localization radius varies among $r_0\in\{1,2,3,4,5,8,15,18\}$. All figures and tables below use an ensemble size of $m=10$, for which
the empirical ensemble covariance has rank at most $m-1=9$, smaller
than the $14$-dimensional unstable--neutral subspace. Tables~\ref{tab:EnKBF_RMSE} and~\ref{tab:EnKBS_RMSE} report the filter and smoother root-mean-square errors (RMSEs), respectively, as functions of the localization radius $r_0$ and the inflation factor $\delta^2$. Here, the RMSE is defined by \[
\mathrm{RMSE}:=\sqrt{\frac{1}{nK}\sum_{k=1}^K\left\|\overline{\vec{x}} (t_{k}) -\vec{x}_{\mathrm{ref}}(t_{k})\right\|^2},
\]
where \(K\) is the total number of time steps and \(\overline{\vec{x}}\) denotes the ensemble mean (filter or smoother). The results show that, with localization and mild inflation, the smoother
consistently improves upon the filter on Lorenz-96. While a suitable
inflation factor slightly reduces the error, localization has a
pronounced impact on the stability and convergence of EnKBS. Compared
with EnKBF, EnKBS is more sensitive to ensemble size. Without
localization or inflation, a stable EnKBS requires an ensemble size of $m\gtrsim31$, whereas the EnKBF
can remain stable with $m=17$. With localization, EnKBS remains stable
with much smaller ensembles, down to $m=5$. These qualitative stability
conclusions are consistent across multiple random-seed tests. Figure~\ref{fig:trj_l96} compares the time trajectories of \(x_1\) between the reference solution and the EnKBF/EnKBS mean, with parameters \(r_0=4\) and \(\delta^2=1.01\).

\begin{table*}[t]
    \caption{Lorenz-96: RMSE of the localized and inflated EnKBF over the evaluation window \([20,100]\). Entries marked as NaN indicate numerical divergence. The symbol $*$ marks the smallest RMSE among the tested parameters.}
    \label{tab:EnKBF_RMSE}
\begin{tabular}{ccccccccc}
    \tophline
         $\delta^2\backslash r_0$&  1&  2& 3& 4& 5& 8&15&18\\
         \middlehline
         1&  0.714&  0.662& 0.656&  0.660&  0.668& 0.704&   0.734& 1.025\\
         \middlehline
         1.0001&   0.712&  0.662& 0.656& 0.659& 0.668& 0.704& 0.733& 1.025\\
         \middlehline
        1.001& 0.700& 0.659& 0.654&  0.658& 0.666& 0.703& 0.732& 1.026\\
        \middlehline
        1.005& 0.683&  0.664& 0.654$*$& 0.654& 0.661& 0.698& 0.728& 1.033\\
        \middlehline
        1.01&  NaN&  0.762& 0.683& 0.667& 0.671& 0.698& 0.727& 1.041\\
        \middlehline
 1.02& NaN& NaN& NaN& NaN& NaN& 0.848& 1.037&1.062\\
\bottomhline
\end{tabular}
\end{table*}

\begin{table*}[t]
    \caption{Lorenz-96: RMSE of the localized and inflated EnKBS over the evaluation window \([20,100]\). Entries marked as NaN indicate numerical divergence. The symbol $*$ marks the smallest RMSE among the tested parameters.}
    \label{tab:EnKBS_RMSE}
\begin{tabular}{ccccccccc}
    \tophline
         $\delta^2\backslash r_0$&  1&  2& 3& 4& 5& 8&15&18\\
         \middlehline
         1&  0.583&  0.567& 0.562&  0.574&   0.591&  0.647&   0.686& NaN\\
         \middlehline
         1.0001&   0.581&  0.560& 0.561&  0.573& 0.590& 0.646& 0.685& NaN\\
         \middlehline
        1.001& 0.565& 0.552& 0.555&  0.568& 0.585& 0.643& 0.682& NaN\\
        \middlehline
        1.005& 0.663&  0.535& 0.531& 0.544& 0.563& 0.630& 0.673& NaN\\
        \middlehline
        1.01&  NaN&  0.643& 0.528& 0.519$*$& 0.541& 0.614& 0.662& NaN\\
        \middlehline
        1.02& NaN& NaN& NaN& NaN& NaN& 0.766& 1.001&NaN\\
    \bottomhline
\end{tabular}
\end{table*}

\subsection{Causal Inference in a Nonlinear Dyad Model}\label{subsec:ACI}

We incorporate the EnKBS into the recently developed assimilative causal inference (ACI) framework \citep{andreou2026assimilative}. Unlike traditional approaches based solely on time-averaged statistics or model forecasts, ACI casts causal inference as an inverse problem through Bayesian data assimilation, enabling time-dependent diagnosis of causal dynamics. In this case study, we use a nonlinear dyad model to demonstrate that the EnKBS-ACI framework is able to infer the causal intensity and temporal influence range of a hidden driver from observations of its effect, dynamically identifying the instantaneous cause-and-effect role interplay.

The key idea behind ACI is that, if a variable \(\vec{y}\) drives another \(\vec{x}\), then conditioning on the effect \(\vec{x}\)'s future evolution should retrospectively improve the estimate of the driver \(\vec{y}\)'s current state.
In ACI, the instantaneous causal influence intensity is quantified by an information-gain metric (ACI metric) based on the relative entropy between the smoother distribution and the real-time filtering distribution.
By further evaluating the information gain across lagged smoothers of different assimilating time windows, we can also obtain the (forward) causal influence range (CIR) of a relationship, which characterizes how far the causal signal propagates in time. This CIR is objective, i.e., it does not rely on subjective or empirical threshold choices.
Here, using the EnKBS to approximate the filter and smoother distributions removes the need to directly observe the driver for causal inference, as long as its effect is observed.

\subsubsection{A Nonlinear Dyad Model with Extreme Events}
We consider a nonlinear dyad model that exhibits intermittent extreme events through a trigger--feedback mechanism:
    \begin{subequations}\label{eq:dyad}
        \begin{align}
            \dot u&=(-d_u+cv)u+F_u+\sigma_u\dot W_u,\label{eq:dyad_u}\\
            \dot v&=-d_vv-cu^2+F_v+\sigma_v\dot  W_v.\label{eq:dyad_v}
        \end{align}
    \end{subequations}
Here, $W_u$ and $W_v$ are independent standard Wiener processes. This physics-constrained low-order model produces extreme events through a trigger--feedback cycle: $v$ modulates the damping of $u$, while $u$ feeds back onto $v$ \citep{chen2020predicting}. Specifically, the growth rate of $u$ is controlled by the time-varying coefficient $-d_u+cv$. When intermittent growth of $v$ exceeds the threshold $d_u/c$, the system enters an anti-damping regime and $u$ grows rapidly. As $u$ amplifies, the term $-cu^2$ suppresses $v$, which restores damping in $u$ and terminates the burst in its amplitude, ultimately conserving the energy of the system.

\subsubsection{The Experiment Setup}
\begin{figure*}[t]
    \includegraphics[width=12cm]{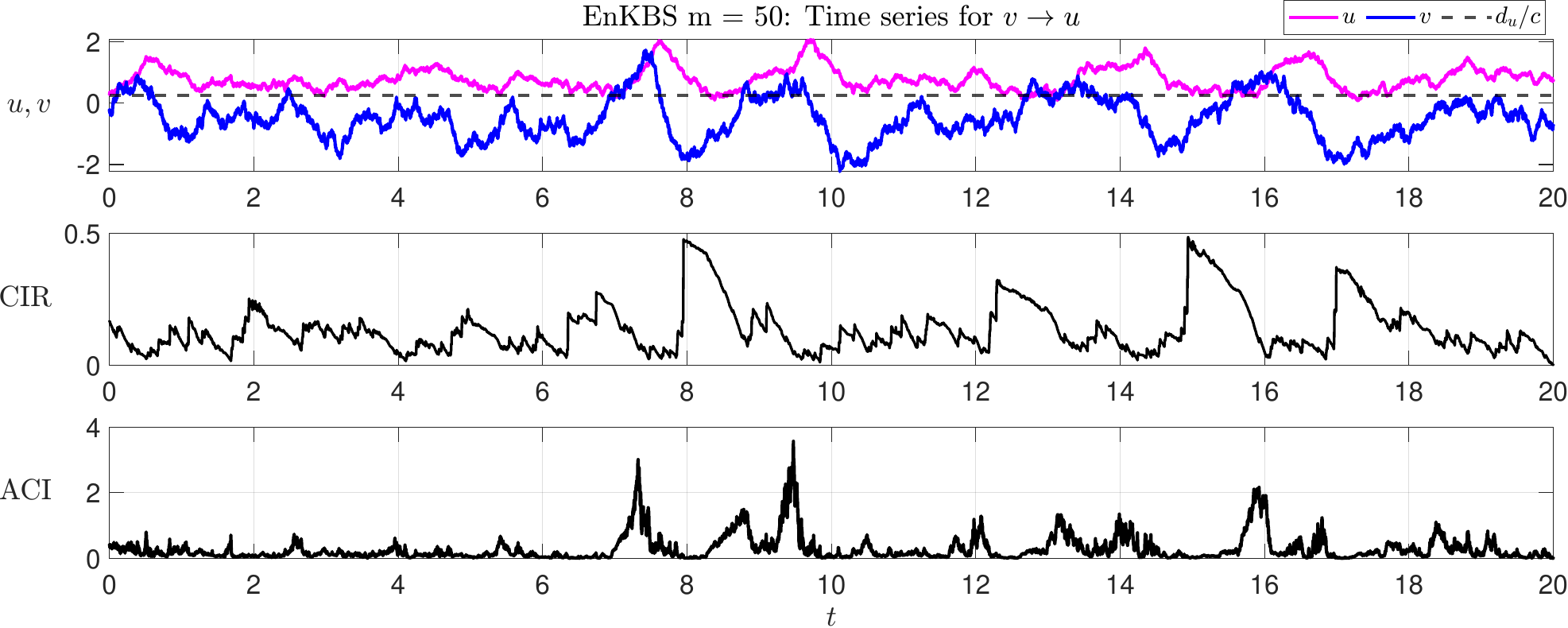}
    \caption{$v\rightarrow u$. CIR and ACI metric computed using the EnKBS with ensemble size $m=50$. Top: true trajectories of $u$ (magenta) and $v$ (blue). The dashed line is the anti-damping threshold \(d_u/c\). Middle: the CIR length for $v\rightarrow u$. Bottom: the ACI metric as a function of time. }
    \label{fig:EnKBSm50ytox}
\end{figure*}

The following parameter values are adopted in Eq.~(\ref{eq:dyad}):
\begin{gather*}
    d_u=0.5,\quad F_u=1,\quad \sigma_u=0.5,\quad c=2,\\
    d_v=0.5,\quad F_v=0.8,\quad \sigma_v=1.
\end{gather*}
We consider both causal directions, the triggering direction $v\to u$ and the feedback direction $u\to v$.
For $v\to u$, the filter and smoother estimate the hidden state \(v\) given the observations of \(u\). Notice that when $u$ is observed, the $v$ dynamics Eq.~(\ref{eq:dyad_v}) become conditionally linear and Gaussian, fitting the conditional Gaussian nonlinear system (CGNS) structure \citep{chen2018conditional,liptser1978statistics}. This provides an analytic filter and smoother that enables a benchmark comparison with the EnKBF and EnKBS \citep{andreou2024martingale,andreou2026assimilative}.
For $u\to v$, we estimate the hidden state $u$ from observations of $v$. Note that this predictability arises through the coupling term $-cu^2$ in Eq.~(\ref{eq:dyad_v}). However, because the map $u\mapsto u^2$ is non-injective, conditioning on $v$ alone does not resolve the sign of $u$. As a result, when the true $u$ is close to $0$, ensemble analysis members may split into opposite branches, artificially inflating the variance and increasing the estimation error. Hence, we fix a reference solution of length $T=30$ over which $u$ does not change sign for both directions. See the top panel of Figure~\ref{fig:EnKBSm50ytox}. We observe that inflation does not noticeably change the results in this regime, and we therefore do not apply inflation.

\subsubsection{Results}
\begin{figure}[t]
    \includegraphics[width=8.3cm]{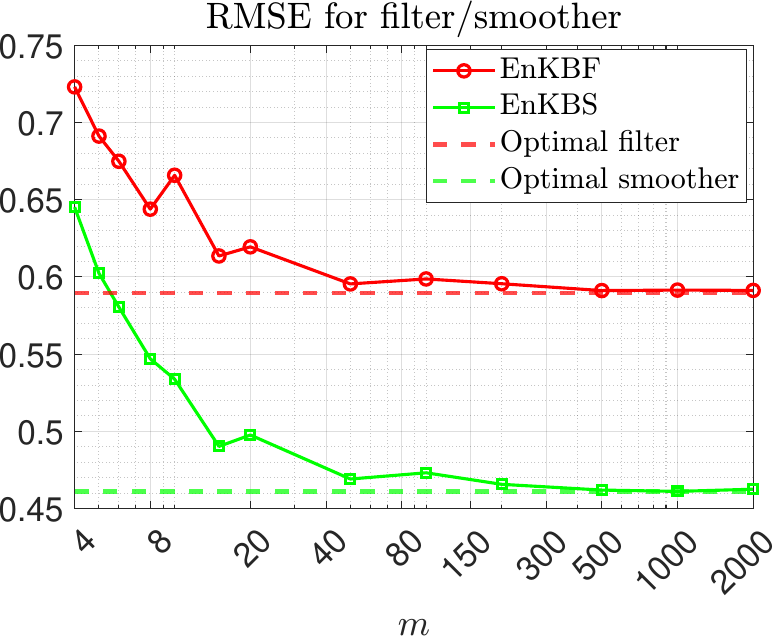}
    \caption{RMSE comparison between EnKBF/EnKBS and CGNS exact filter/smoother as a function of ensemble size \(m\). Red circle: EnKBF. Green square: EnKBS. Red dashed: optimal filter. Green dashed: optimal smoother. The RMSE between the filter/smoother estimates and the truth of \(v\) is evaluated over a time window of \([0,500]\).}
    \label{fig:RMSE}
\end{figure}

\begin{figure*}[t]
    \includegraphics[width=12cm]{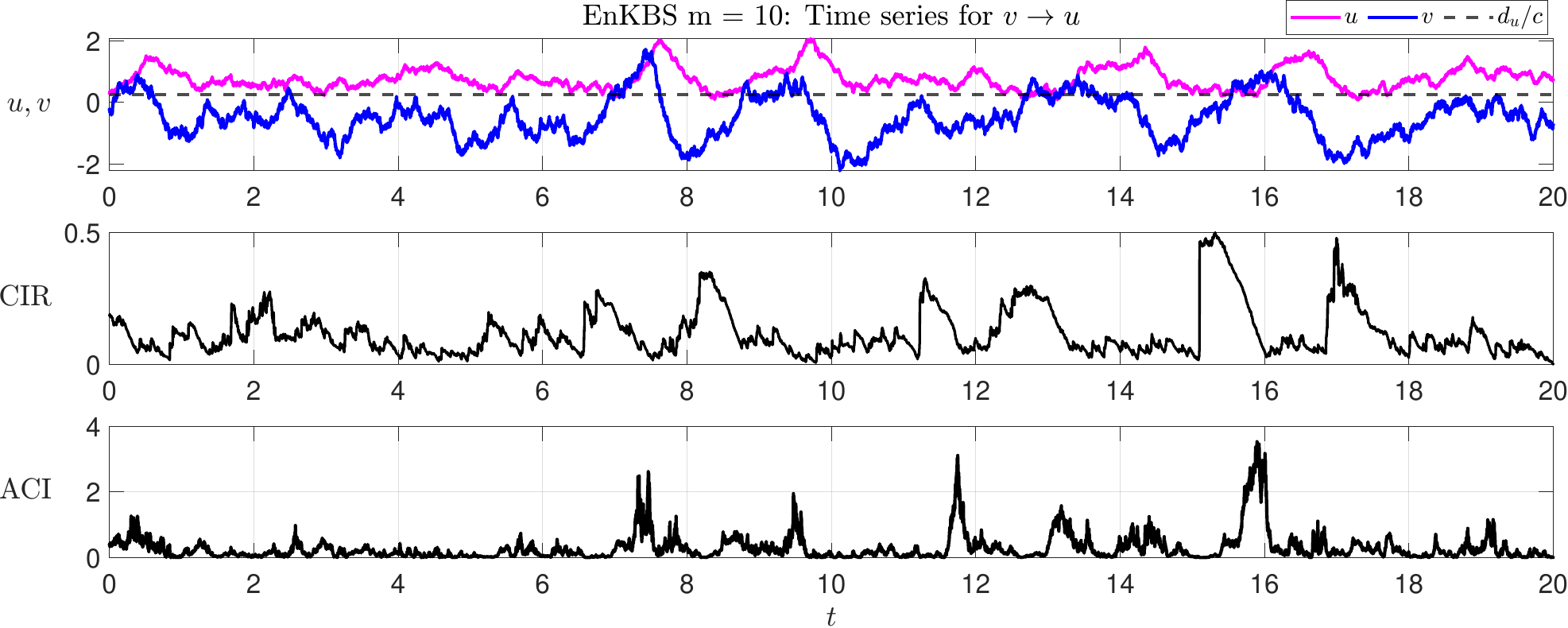}
    \caption{$v\rightarrow u.$ CIR and ACI metric computed using EnKBS with ensemble size $m=10$. Top: true trajectories of $u$ (magenta) and $v$ (blue). The dashed line is the anti-damping threshold \(d_u/c\). Middle: the CIR length for $v\rightarrow u$. Bottom: the ACI metric as a function of time.}
    \label{fig:EnKBSm10ytox}
\end{figure*}

\noindent\textbf{Triggering: \( {v\rightarrow u}\).} Figure~\ref{fig:RMSE} compares the RMSE between the EnKBF/EnKBS and CGNS analytic filter/smoother estimates of \(v\) as a function of ensemble size \(m\).
As \(m\) increases, the EnKBS ensemble mean and covariance converge to the analytic statistics, with the resulting ACI diagnostics converging accordingly. Figures~\ref{fig:EnKBSm50ytox} and \ref{fig:EnKBSm10ytox} show the ACI results with an ensemble size of \(m=50\) and \(m=10\) respectively.  With $m=50$, the EnKBS diagnostics closely match the analytic benchmark. We therefore display only the CIR and ACI metric results with the EnKBS.  For \(m=10\), the ACI metric does not preserve relative peak magnitudes and underestimates the CIR at the beginning of the rise of $v$, i.e., during the build-up phase of $u$'s extreme event; this is especially pronounced at prolonged or delayed onsets. Nevertheless, it captures the peak locations and yields qualitatively similar patterns for the CIR and ACI metric.

In general, the significant ACI metric coincides with significant $v$ before extreme events of \(u\), e.g., $t\approx 7$, $9.5$, and $16$, indicating the strongest instantaneous influence where \(v\) acts as the primary anti-damping driver of these extreme events. The CIR peaks at the onset of $v$'s growth, e.g., $t\approx 6.5$, $8$, and $15$, suggesting long-range dependence from \(v\) to \(u\) during the onset stage of these extreme events, and the triggering conditions of extreme events of \(u\) are established well in advance. The CIR then shortens, and the ACI metric drops sharply as the events approach their peaks and begin to decay.

From a data assimilation perspective, the filter cannot detect the triggering conditions at the onset phases unless they are manifested in the observations. In contrast, the smoother is able to clearly improve upon the filter with future information in hindsight. As a result, a larger ACI metric and a longer CIR are obtained. Later, when \(u\) begins to strengthen, near-future observations already provide sufficient information for effective smoothing. Once \(u\) becomes strong enough, its signal-to-noise ratio enables the filter to effectively capture the relevant triggering dynamics, leaving little room for the smoother to improve estimates. Therefore, the large CIR and ACI metric does not persist in the peak and demise stages of extreme events of \(u\).

\begin{figure*}[t]
    \includegraphics[width=12cm]{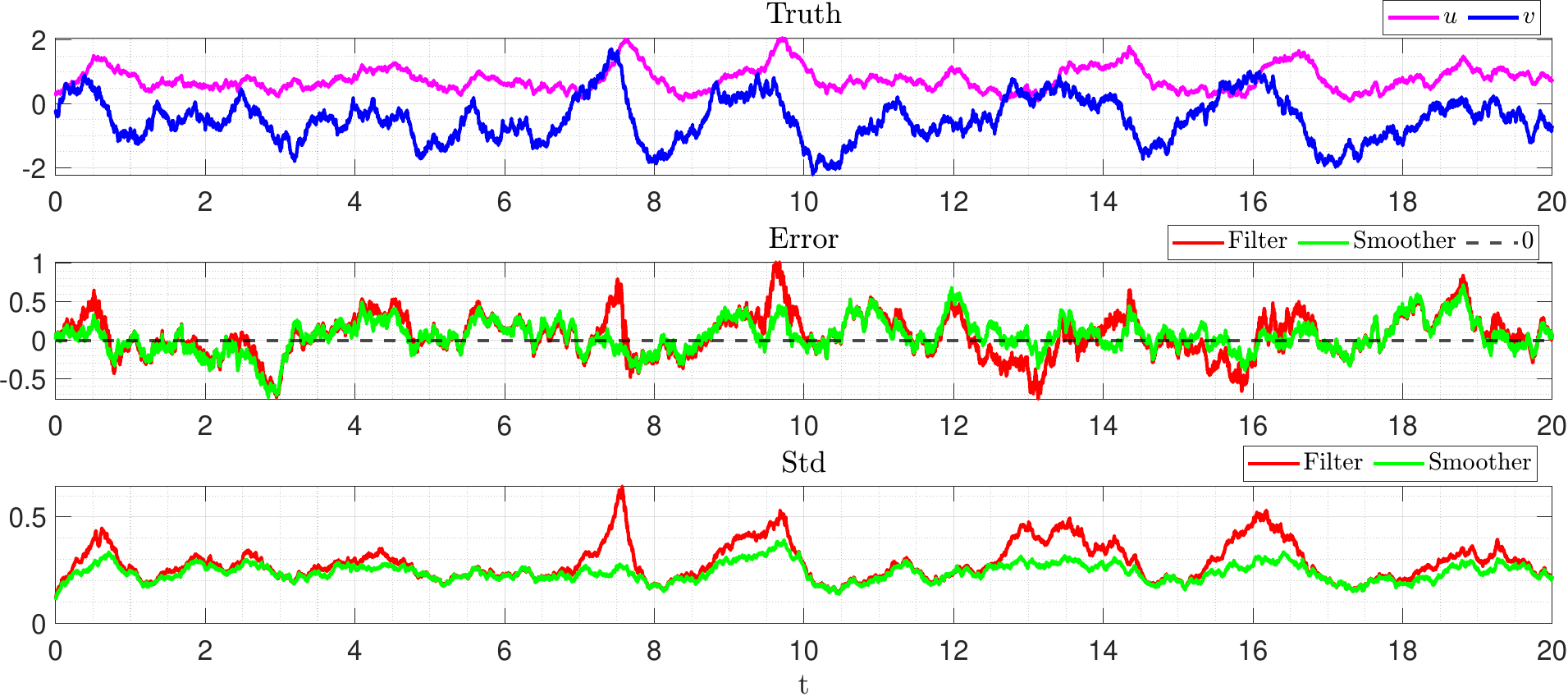}
    \caption{$u\rightarrow v$. EnKBF vs. EnKBS diagnostics with ensemble size $m=50$. Top: true trajectories of $u$ (magenta) and $v$ (blue). Middle: estimation errors of $u$ for the filter (red) and smoother (green), computed as $u_{\mathrm{ref}}-u_{\mathrm{f}}$ and $u_{\mathrm{ref}}-u_{\mathrm{s}}$, respectively. Values closer to \(0\) (dashed) indicate smaller errors. Bottom: empirical standard deviation of $u$ for the filter (red) and smoother (green).}
    \label{fig:good_x_to_y_err_m50}
\end{figure*}

\noindent\textbf{Feedback: \( {u\rightarrow v}\).}
In this regime, experiments use an ensemble size of \(m=50\) (based on the results of Figure~\ref{fig:RMSE}). Figure~\ref{fig:good_x_to_y_err_m50} shows the estimation errors and empirical standard deviations for both the filter and the smoother. During extreme events, the smoother provides a clear correction relative to the filter. No sign-splitting is observed since $u>0$ throughout.

\begin{figure*}[t]
    \includegraphics[width=12cm]{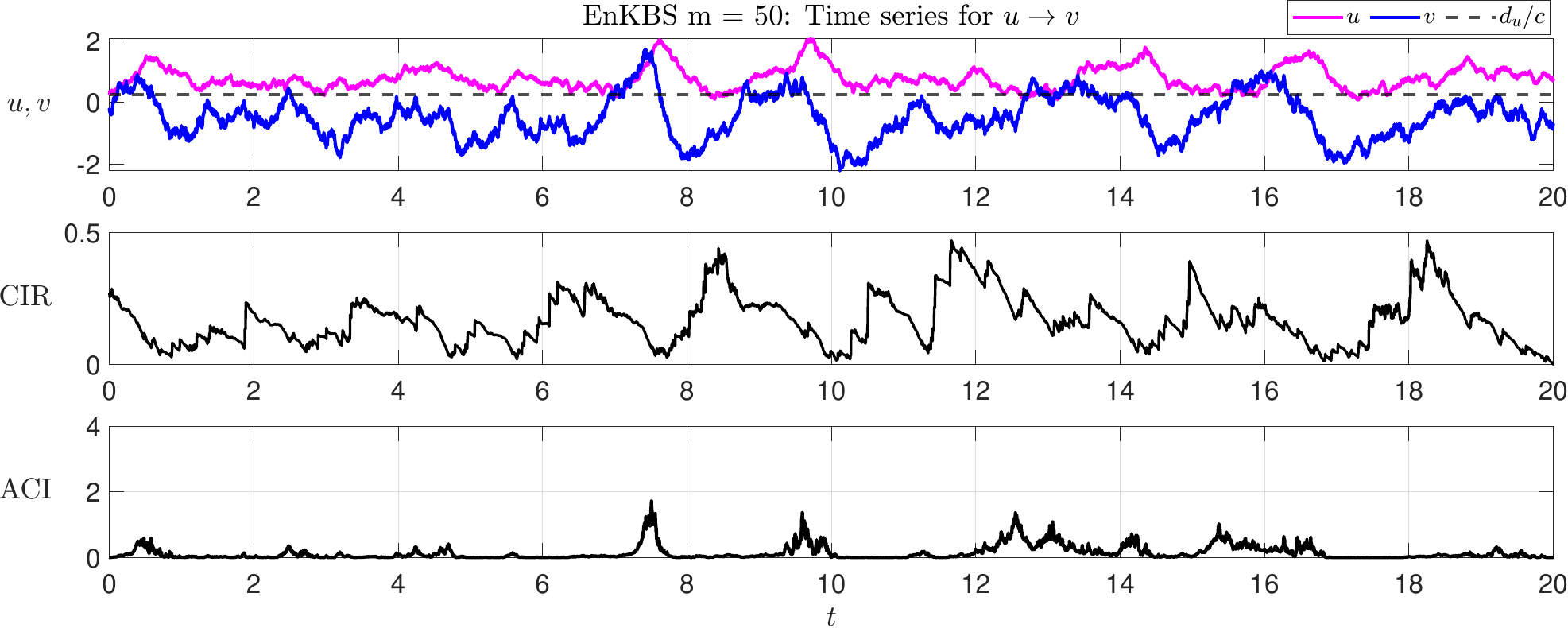}
    \caption{$u\rightarrow v$. CIR and ACI metric computed using EnKBS with ensemble size $m=50$. Top: true trajectories of $u$ (magenta) and $v$ (blue). The dashed line is the anti-damping threshold \(d_u/c\). Middle: the CIR length for $u\rightarrow v$. Bottom: the ACI metric as a function of time.}
    \label{fig:good_x_to_y_ACI_m50}
\end{figure*}

Figure~\ref{fig:good_x_to_y_ACI_m50} shows the CIR and ACI metric. Large ACI metric matches the peak of $v$, e.g., $t\approx 7.6$, $9.6$, and $13$. This indicates the strongest causal direction in which $u$ suppresses the growth of $v$. The CIR is relatively long during the onset stages of extreme events of \(u\), e.g., $t\approx 6.5$, $8.5$, $10.8$, and $12$, suggesting that, as in the \(v\to u\) direction, the \(u\to v\) feedback has been established well before extreme events. As $v$ approaches its peak and proceeds towards its demise, the CIR gradually shortens. Changing the ensemble size to $m=100$ or $m=10$ yields quantitatively similar results (not shown).

To interpret these patterns, given a time series of $\mathrm{d} v$ in the $v$-equation Eq.~(\ref{eq:dyad_v}), the filter infers $u$ through the negative coupling term $-cu^2$. As a result, the filter cannot anticipate the growth and peak of $u$ unless the subsequent decline of $v$ is observed. In contrast, the smoother can trace back from this decline to correctly estimate \(u\). Consequently, both the CIR and the ACI metric become large during the onset stage, and the CIR gradually shortens as the decline of \(v\) is imminent.  Later, $v$ starts to decay. This phenomenon already provides an informative signal to the filter, and the information flow $u\rightarrow v$ has largely been realized and manifested in the CIR and ACI metric during the onset period. Therefore, the filter can achieve an accurate estimate via dynamical interpolation, so the additional future information becomes negligible. As a result, the ACI metric does not persist after the peak and through the full demise of the extreme events.

\subsection{Data-Driven Discovery of Dynamics with EnKBS-Based Sampling}\label{subsec:est}

In this subsection, we adapt a causality-based learning algorithm to recover the underlying dynamics of partially observed systems.
The key idea is to formulate model learning as an iterative procedure that alternates among conditional sampling, structure identification, and parameter estimation \citep{chen2023causality}. The procedure consists of three steps:
\begin{enumerate}
  \item \textbf{Conditional sampling.} Based on partially observed time series, reconstruct trajectories of the hidden variables.
  \item \textbf{Structure identification.} Detect which candidate terms are active in the governing equations by quantifying the information transfer from each functional, in a library of functions, to the trajectory of each state.
  \item \textbf{Parameter estimation.} Given the identified structure, recover the corresponding coefficients via maximum likelihood estimation.
\end{enumerate}

The learning problem considered here is offline: the complete
observation trajectory \(\vec{y}_{[0,T]}\) is available before each
learning iteration. Smoothing is therefore a natural choice for the
conditional-sampling step, allowing both past and future observations
to reconstruct the hidden trajectory.

The EnKBS plays a key role in Step~1.
In our approach, the smoother ensemble member trajectories provide conditional samples of the hidden states. We treat these trajectories as synthetic observations to complete the dataset and enable subsequent identification of the model structure and parameters.

We incorporate physics constraints, e.g., energy conservation in quadratic terms, to prevent nonphysical finite-time blow-up and stabilize learning \citep{majda2012physics,harlim2014ensemble}. These constraints are enforced in the parameter-estimation step together with other constraints to facilitate learning convergence.

\subsubsection{The Lorenz-84 System}

Numerical tests are performed on the Lorenz-84 model \citep{lorenz1984irregularity}:
    \begin{equation}\label{eq:l84_truth}
        \begin{aligned}
            \dot  x &= -(y^2+z^2) - a(x-f) + \sigma_x \dot W_x,\\
           \dot y&= -b x z + x y - y + g + \sigma_y \dot W_y,\\
            \dot z &= b x y + x z - z + \sigma_z \dot W_z.
        \end{aligned}
    \end{equation}
We use the following parameter setting to induce chaotic dynamics:
    \begin{equation*}
        a=\frac14,\qquad b=4,\qquad f=8,\qquad g=1,\qquad \sigma_x=\sigma_y=\sigma_z=0.1.
    \end{equation*}
Equation~(\ref{eq:l84_truth}) is a low-order atmospheric conceptual model: $x$ denotes the intensity of the mid-latitude westerly zonal flow, while $y$ and $z$ represent the cosine and sine longitudinal phases of superimposed large-scale wave patterns.
Stochastic forcing is included to represent interactions with unresolved sub-grid processes, while $f$ and $g$ represent symmetric (e.g., by solar heating) and asymmetric (e.g., from topographic influences) cross-latitude heating contrasts, respectively.

\subsubsection{The Experiment Setup}
We observe $(y,z)$ and treat $x$ as hidden.
We use the following 12-term candidate function library:
    \begin{equation*}
1,\;y,\;z,\;y^2,\;z^2,\;yz,\;x,\;xy,\;xz,\;xy^2,\;xz^2,\;xyz.
    \end{equation*}
This library matches the CGNS structure, enabling direct comparison with the CGNS optimal sampling benchmark \citep{chen2023causality,chen2018conditional}.
In the structure-identification step, the constant term is always included by default.

We use the following initial model guess, as in \citep{chen2023causality}:
    \begin{equation*}
        \begin{aligned}
           \dot  x &= y^2 - z^2 + 2 + (y^2-z^2)x + \sigma_x \dot W_x,\\
           \dot  y &= -y - 2y^2 + z^2 + 1 + (-y-8z-yz)x + \sigma_y \dot W_y,\\
          \dot  z   &= -z + z^2 - yz + (8y+z+z^2)x + \sigma_z \dot W_z.
        \end{aligned}
    \end{equation*}
The initial guesses for the noise amplitudes of observed variables are set to $\sigma_y=\sigma_z=1$, while the noise of the hidden \(x\) is treated as known and fixed at  $\sigma_x=0.1$ for simplicity.
The time window length is $T=500$ with numerical time step $\Delta t=0.001$. We run $120$ learning iterations. For EnKBS-based sampling, we use an ensemble size of $m=50$. Using $m=100$ or $m=10$ produces similar learning behavior, except that the $m=10$ case typically converges a few iterations later.

\subsubsection{Results}
\begin{figure*}[t]
    \includegraphics[width=12cm]{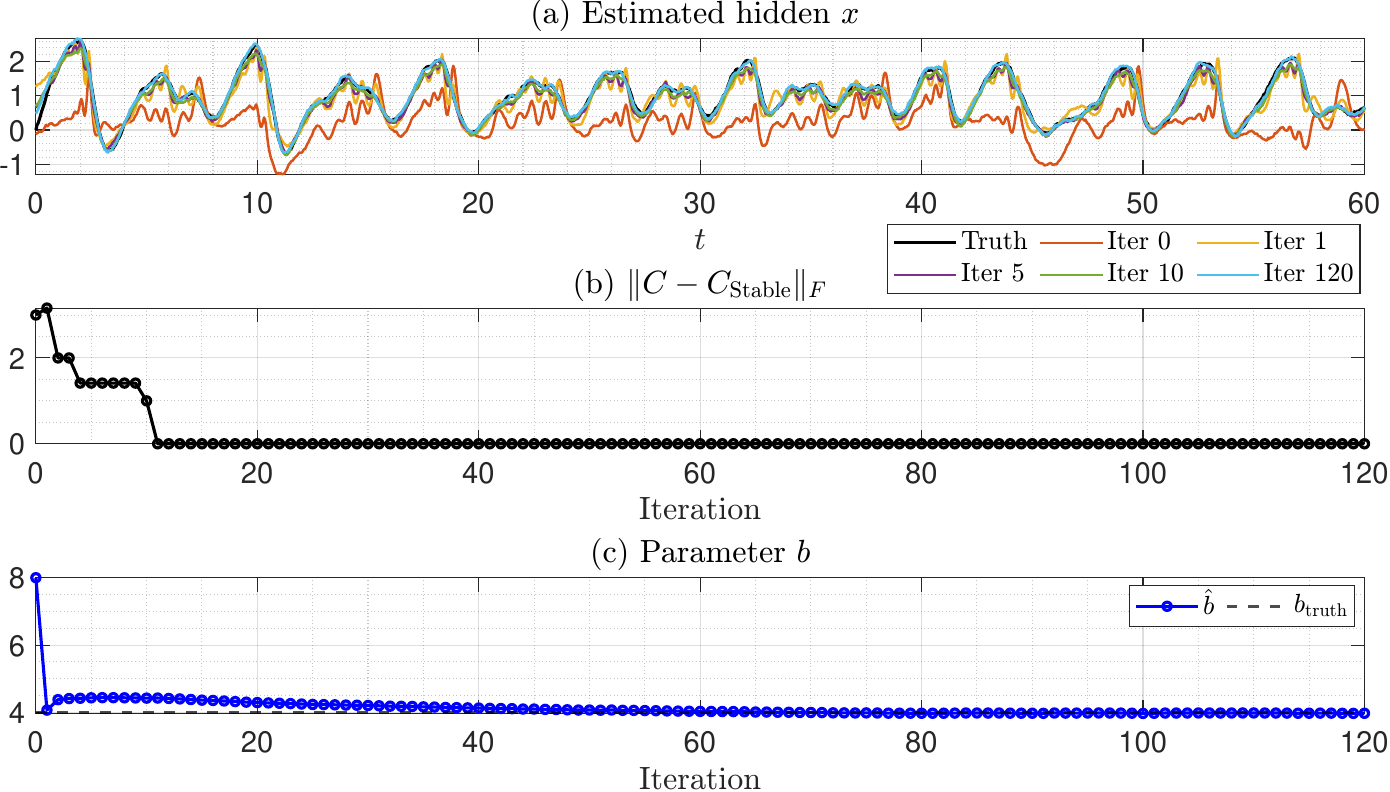}
    \caption{Learning progress for Lorenz-84 using EnKBS-based sampling. (a) Sampled trajectories of the hidden state $x$ at selected iterations compared with the truth; iteration~0 corresponds to the initial guess. (b) Structure error measured by Frobenius norm $\|\mathbf C-\mathbf C_{\mathrm{stable}}\|_F$. (c) Evolution of the estimated parameter $\hat b$ toward the true value $b$.}
    \label{fig:newbaseest}
\end{figure*}

Figure~\ref{fig:newbaseest} shows that, with EnKBS-based sampling in Step~1, the iterative procedure identifies the correct model structure, converges, and recovers the hidden trajectory $x(t)$.
In particular, Figure~\ref{fig:newbaseest}(a) shows that the sampled trajectory of the unobserved variable $x$ gradually approaches the truth as iterations proceed.

To quantify structure convergence, we introduce a causation entropy indicator matrix $\mathbf C\in\{0,1\}^{N\times M}$ with system dimension $N=3$ and  the number of candidate functions $M$ \citep{chen2023causality}.
Each entry is obtained by thresholding the causation entropy:
    \begin{equation*}
        \mathbf C_{i,j}=
        \begin{cases}
            1, & \mathrm{CE}_{i,j}\ge r,\\
            0, & \mathrm{CE}_{i,j}< r,
        \end{cases}
        \qquad r=10^{-3},
    \end{equation*}
where $\mathrm{CE}_{i,j}$ denotes the causation entropy from the $j$-th candidate function to the $i$-th state equation \citep{elinger2021causation,sun2014causation}.
We then compute the Frobenius norm of the matrix difference $\|\mathbf C-\mathbf C_{\mathrm{stable}}\|_F$, where $\mathbf C_{\mathrm{stable}}$ denotes the indicator matrix after the identified structure has stabilized to the true model structure Eq.~(\ref{eq:l84_truth}).
The decay in Figure~\ref{fig:newbaseest}(b) indicates that the model structure is fully identified after about 10 iterations.

Figure~\ref{fig:newbaseest}(c) shows that the estimate of $\hat b$ converges to the true value $b$ after roughly 60 iterations.
After that point, the remaining coefficients show no noticeable changes. The final parameter estimates at the 120th iteration are listed in Table~\ref{tab:l84_table1}. The identified coefficients are close to the true values.

\begin{table}[t]
    \caption{Coefficients in the true Lorenz-84 system Eq.~(\ref{eq:l84_truth}) and the identified model at the 120th iteration. The identified model has the same structure as the truth. Here, $\theta^x_{xy}$ denotes the coefficient of the feature $xy$ in the $x$-equation, and analogous notation is used for other terms.}
    \label{tab:l84_table1}
\begin{tabular}{lcccccc}
\tophline
 & $\theta^x_{x}$ & $\theta^y_{y}$ & $\theta^z_{z}$ & $\theta^x_{y^2}$ & $\theta^x_{z^2}$ & $\theta^y_{xz}$\\
 \middlehline
Truth      & -0.25& -1& -1& -1& -1& -4\\
\middlehline
Identified & -0.2459& -1.0135& -1.0037& -1.0023& -1.0011& -3.9732\\
\middlehline
 & $\theta^y_{xy}$ & $\theta^z_{xy}$ & $\theta^z_{xz}$ & $\theta^x_{1}$ & $\theta^y_{1}$ & $\theta^z_{1}$ \\
 \middlehline
Truth      & 1& 4& 1& 2& 1& 0\\
\middlehline
Identified & 1.0023& 3.9732& 1.0011& 1.9946& 1.0088& -0.0061\\
\bottomhline
\end{tabular}
\end{table}

\begin{figure*}[t]
    \includegraphics[width=12cm]{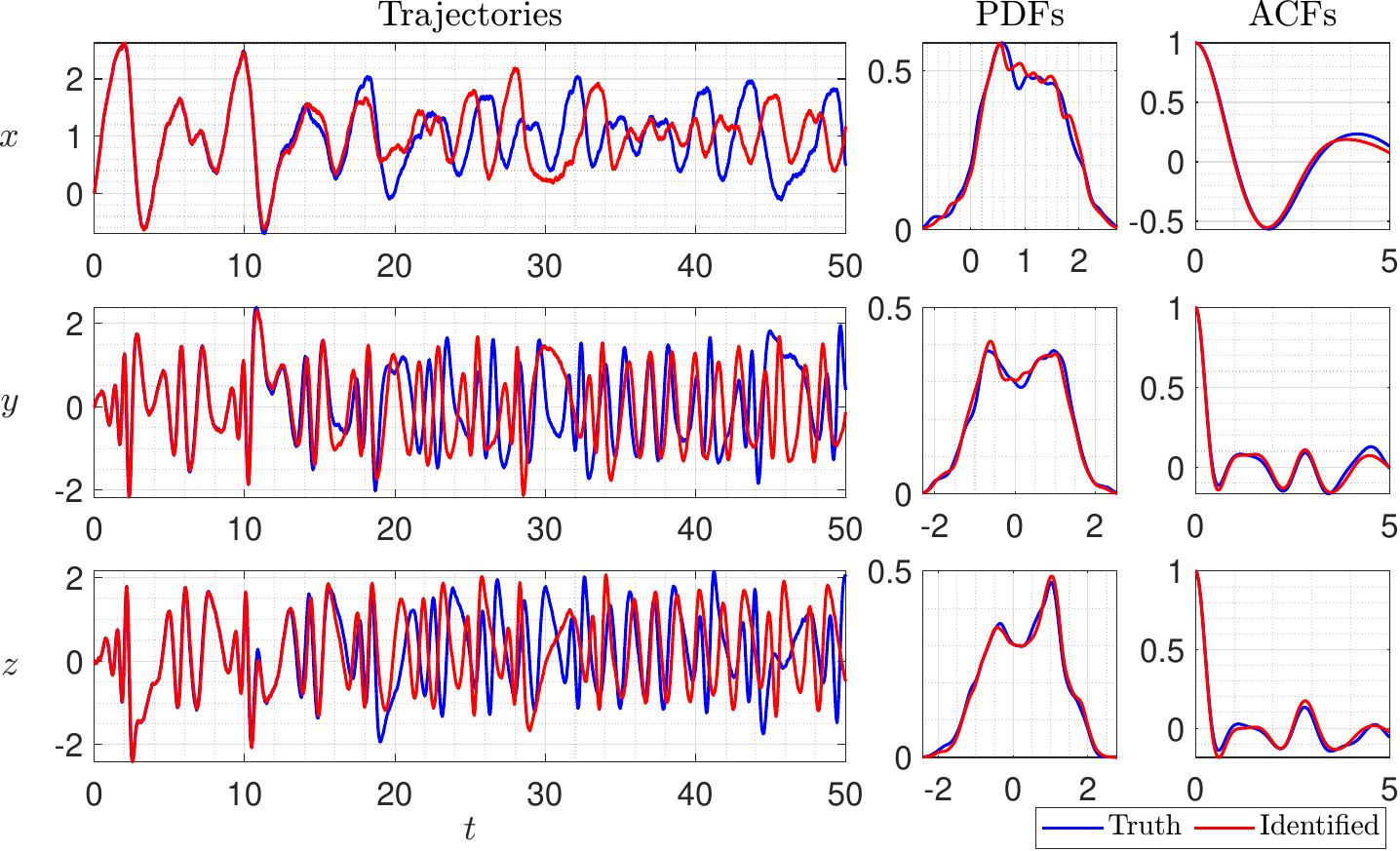}
    \caption{EnKBS-based sampling. Trajectory and statistics comparison between the true Lorenz-84 system (blue) and the identified model (red). In addition to short-term path-wise agreement, the identified model reproduces the PDF and the temporal ACF of the truth.}
    \label{fig:newbasetraj}
\end{figure*}

Figure~\ref{fig:newbasetraj} compares trajectories and two key statistics generated by the true system and the identified model. When generating the time series, we use the same random seed and initial conditions for the true system and the identified model. The two trajectories match closely up to about $t\approx 10.5$, demonstrating short-term path-wise consistency. Due to chaos, we do not expect long-term point-wise agreement over the full window.
Moreover, the identified model reproduces both the probability density function (PDF) and the temporal autocorrelation function (ACF), indicating statistical agreement with the truth.

\begin{figure*}[t]
    \includegraphics[width=12cm]{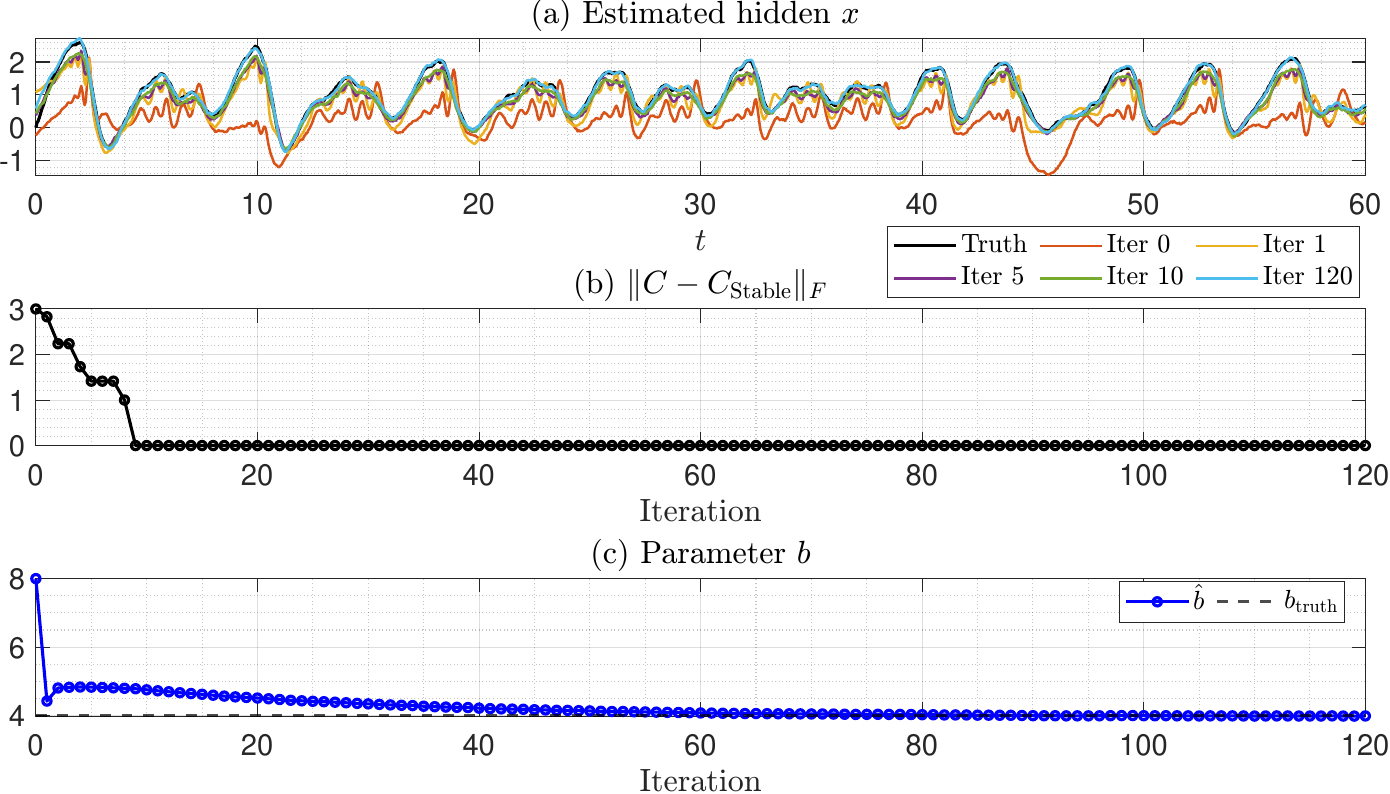}
    \caption{Learning progress for Lorenz-84 using CGNS optimal sampling (benchmark). (a) Sampled trajectories of the hidden state $x$ at selected iterations compared with the truth; iteration~0 corresponds to the initial guess. (b) Structure error measured by Frobenius norm $\|\mathbf C-\mathbf C_{\mathrm{stable}}\|_F$. (c) Evolution of the estimated parameter $\hat b$ toward the true value $b$.}
    \label{fig:CGest}
\end{figure*}

\begin{figure*}[t]
    \includegraphics[width=12cm]{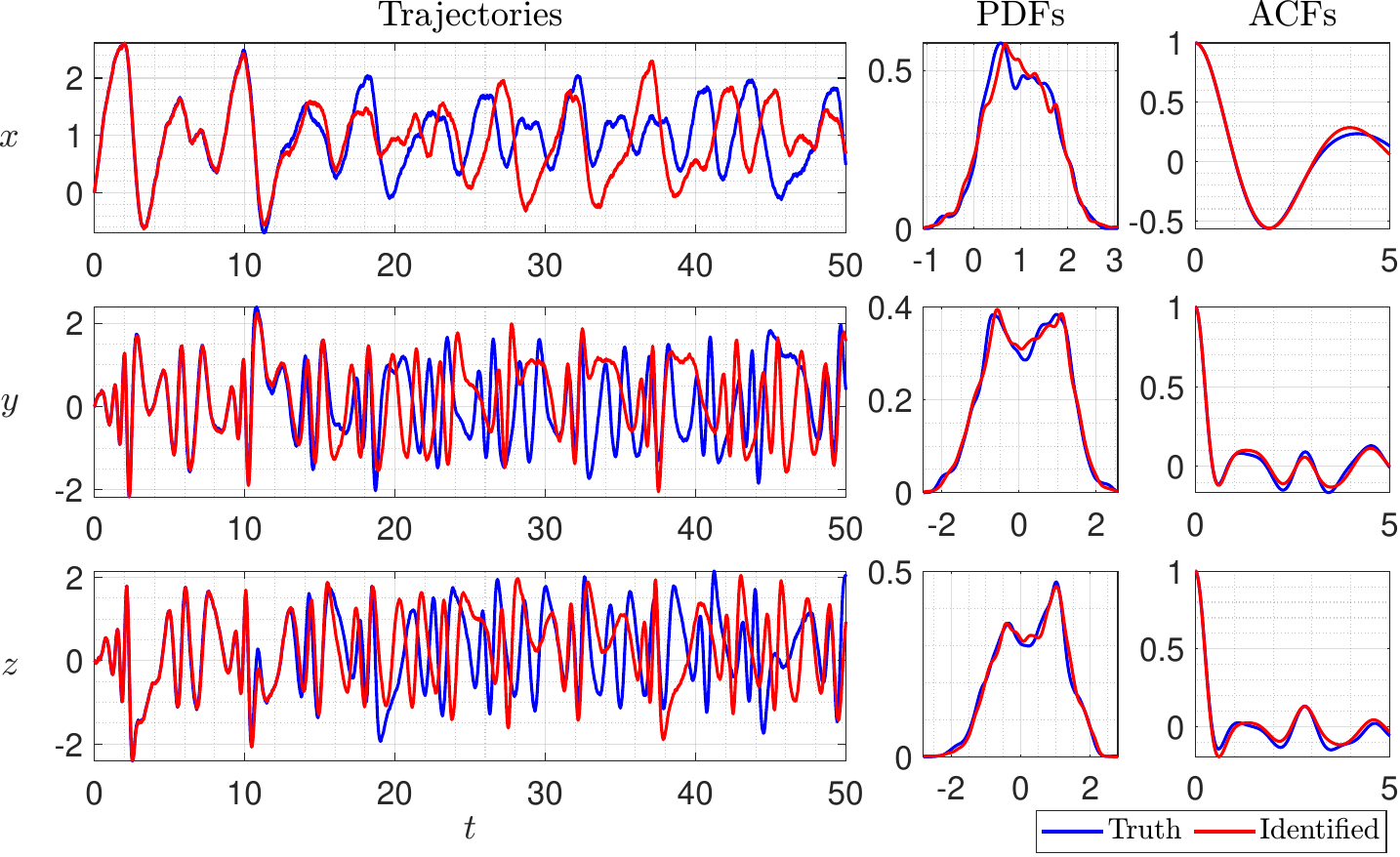}
    \caption{CGNS optimal sampling. Trajectory and statistics comparison between the true Lorenz-84 system (blue) and the identified model (red). In addition to short-term path-wise agreement, the identified model reproduces the PDF and the temporal ACF of the truth.}
    \label{fig:CGtraj}
\end{figure*}

Finally, Figures~\ref{fig:CGest} and \ref{fig:CGtraj} show the benchmark results using CGNS optimal sampling. The overall learning behavior is similar to the EnKBS-based results in Figures~\ref{fig:newbaseest} and \ref{fig:newbasetraj}.

\conclusions \label{sec:Con} 
In this paper, we formulate a systematic ensemble Kalman--Bucy smoother (EnKBS) framework for state estimation in nonlinear dynamical systems.
The proposed method consists of a forward--backward architecture: a forward filtering pass for online estimation, followed by a backward smoothing pass for retrospective refinement.
To handle general nonlinear systems, we adopt an ensemble-based approximation in which conditional moments are replaced by empirical ensemble moments, avoiding tangent-linear and adjoint models.
The forward pass assimilates observations to produce filtering estimates, while the backward pass treats the filtered trajectory as a prior and incorporates future information to refine the state estimates.

Through numerical experiments, we show that the EnKBS consistently improves the filtered estimates by reducing error and uncertainty, converging asymptotically to the optimal smoother when an analytic benchmark is available.
We further demonstrate that covariance localization is important for
numerical stability in high-dimensional settings, while mild inflation
can provide an additional reduction in estimation error.
With regularization, spurious long-range correlations and rank deficiency are mitigated, preventing numerical divergence in high dimensions.
Beyond state estimation, we integrate EnKBS into a Bayesian causal inference framework and a causality-based iterative learning algorithm for partially observed systems \citep{andreou2026assimilative,chen2023causality}.
With an ensemble size of $O(10)$, EnKBS provides effective state estimates and trajectory samples of unobserved variables in nonlinear systems without analytic smoothers, enabling causality diagnostics around extreme events and supporting model discovery.

As a final remark, although measurements are typically available at discrete times, the continuous-time EnKBS provides a complementary perspective to discrete-time and mixed-time smoothers \citep{kalman1961new,bergemann2012Ensemble}.
This viewpoint is consistent with the continuous-time nature of the underlying physics, and is particularly natural when the goal is to discover dynamical mechanisms.
We interpret observed variables as components of the same continuous-time system, rather than as external analysis updates.
Assimilation then exploits intrinsic coupling and information transfer within the system, using observed trajectories to constrain the inferred hidden states.

Future work is twofold: expanding the EnKBS methodology and its dynamical discovery applications.
On the methodology side, while EnKBS performs well on test systems in this paper, these examples largely involve quadratic or conditionally quadratic structures.
A natural next step is to evaluate EnKBS on more realistic, high-dimensional, and strongly nonlinear models.
In particular, we can study accuracy and convergence as well as filter/smoother degeneracy together with possible mitigation strategies.
Another direction is to couple the EnKBS with modern
machine-learning-enhanced data assimilation approaches. When repeated
integration of the dynamical model is computationally expensive, a
learned emulator or data-driven surrogate could be used to propagate
the forecast ensemble while retaining the EnKBF/EnKBS assimilation
framework \citep{hodyss2026diffusion}. Another possibility is to
learn a surrogate for the EnKBS backward correction
\(\mathbf\Sigma\mathbf P_{\mathrm f}^{-1}
(\vec{x}_{\mathrm s}^{(i)}-\vec{x}_{\mathrm f}^{(i)})\). Approximating this term could potentially reduce the cost associated with repeated
construction and inversion of the empirical covariance. This direction is part of our ongoing work on further developing the EnKBS framework. 
On the application side, further experiments and analyses are needed to clarify the mechanisms behind causal information flow and causal reversals in driver-feedback cycles.
For example, alongside the forward CIR, we can investigate with the EnKBS the backward CIR to assess the causal precursors that influenced the currently observed event, enabling attribution of significant events in general complex systems \citep{andreou2025bridging}.


\codeavailability{All MATLAB code used to generate the numerical results and figures reported in this paper is available at \url{https://github.com/jiangzh67/EnKBS} under the MIT License.} 





\appendix
\section{Formal Connection to Discrete-time Formulations}\label{app:der}

In this appendix, we show how the continuous-time stochastic EnKBF formulation Eq.~(\ref{eq:enkbf-sde-stc}) and the EnKBS Eq.~(\ref{eq:enkbs_member_sde}) connect to the standard discrete-time Kalman filter and the RTS smoother \citep{kalman1960new,rauch1965maximum}.
The following formal small-step derivation is carried out for
linear--Gaussian dynamics in the mean-field limit \(m\rightarrow\infty\),
where the empirical ensemble mean and covariance are identified with their
corresponding filtering or smoothing moments \citep{calvello2025ensemble}. 
The nonlinear formulations are then obtained by restoring the
nonlinear model drift and using the ensemble approximations
Eqs.~(\ref{eq:abb_tl})--(\ref{eq:abb_aj}). As a nonrigorous argument, this highlights the limiting equation structure rather than
establishing convergence for general nonlinear systems.

Fix a small observation time step $\tau>0$ and define $t_k=k\tau$.
We use the subscript $(\cdot)_k$ to denote evaluation at time $t_k$, e.g., $\vec{x}_k=\vec{x}(t_k)$.
We consider a first-order semi-implicit Euler--Maruyama type discretization of the continuous model Eq.~(\ref{eq:general_system}) \citep{law2015data}:
\begin{subequations}\label{eq:EM-dis}
        \begin{align}
        \vec{x}_{k+1} &=\vec{x}_k + \tau \vec{f}_k + \sqrt{\tau}\mathbf\Sigma_k^{1/2}\vec{B}_k,\\
        \Delta \vec{y}_{k+1} & = \tau \vec{h}_{k+1} + \sqrt{\tau}\mathbf\Gamma_{k+1}^{1/2}\vec{W}_{k+1},
    \end{align}
\end{subequations}
where $\vec{B}_k$ and $\vec{W}_{k+1}$ are independent standard Gaussian vectors, and $\Delta \vec{y}_{k+1}=\vec{y}_{k+1}-\vec{y}_k$.

For the linear drifts, we have
\begin{equation}\label{eq:tgt_appr}
    \vec f_k=\mathbf F_k\vec x_k,\qquad
    \vec h_{k+1}=\mathbf H_{k+1}\vec x_{k+1}.
\end{equation}
Substituting Eq.~(\ref{eq:tgt_appr}) into Eq.~(\ref{eq:EM-dis}) gives
\begin{subequations}\label{eq:EM-dis_appro}
\begin{align}
    \vec x_{k+1}
    &=(\mathbf I+\tau\mathbf F_k)\vec x_k
      +\sqrt{\tau}\mathbf\Sigma_k^{1/2}\vec B_k,\\
    \Delta\vec y_{k+1}
    &=\tau\mathbf H_{k+1}\vec x_{k+1}
      +\sqrt{\tau}\mathbf\Gamma_{k+1}^{1/2}\vec W_{k+1}.
\end{align}
\end{subequations}

We use the notation $(\cdot)_{k\mid \ell}$ for quantities at time $t_k$ conditioned on observations up to time $t_\ell$.
In this notation, $\vec{x}_{k+1\mid k}$ is the prior obtained by dynamical forecast from $t_k$ to $t_{k+1}$, and $\vec{x}_{k+1\mid k+1}$ is the posterior after assimilating the observation at $t_{k+1}$. Eq.~(\ref{eq:EM-dis_appro}) leads to the following ensemble Kalman filter analysis update for $i=1,\dots,m$ \citep{burgers1998analysis,reich2011dynamical}:
    \begin{equation}\label{eq:filter_correction}
        \vec{x}_{k+1\mid k+1}^{(i)} = \vec{x}_{k+1\mid k}^{(i)} + \mathbf K_{k+1}
        \left(
        \Delta  \vec{y}_{k+1} - \tau \mathbf H_{k+1}\vec{x}_{k+1\mid k}^{(i)}-\sqrt{\tau}\mathbf\Gamma^{1/2}_{{k+1}} \vec{W}^{(i)}_{{k+1}} 
        \right),
    \end{equation}
with \( \sqrt{\tau}\vec{W}^{(i)}_{{k+1}} \sim\mathcal{N}(\mathbf 0, \tau\mathbf{I})\) a discrete Brownian increment and \(\mathbf K_{k+1}\) the ensemble-based Kalman gain:
\[
    \mathbf K_{k+1}:=\mathbf P_{k+1\mid k}\left(\mathbf H_{k+1}\right)^\top\left(\tau\mathbf H_{k+1}\mathbf P_{k+1\mid k}(\mathbf H_{k+1})^\top + \mathbf\Gamma_{{k+1}}\right)^{-1}.
\]
Here, \(\mathbf P_{k+1\mid k}\) is the empirical covariance of the ensemble prior \(\vec{x}_{k+1|k}^{(i)}\). From the analysis step of the Kalman filter, we have \citep{asch2016data,kalman1960new}
\[
    \mathbf P_{ k+1|k+1} = \mathbf P_{k+1\mid k} - \tau \mathbf K_{k+1}\mathbf H_{k+1} \mathbf P_{k+1\mid k}=\mathbf P_{k+1\mid k} +O( \tau).
\]
Together with
\[
    \mathbf K_{k+1}
    =
    \mathbf P_{k+1\mid k+1}\mathbf H_{k+1}^\top
    \mathbf\Gamma_{k+1}^{-1}
    + O(\tau),
\]
Eq.~(\ref{eq:filter_correction}) gives, in the formal continuous-time limit,
\begin{equation}\label{eq:enkf-analysis-limit}
    \vec{x}_{k+1\mid k+1}^{(i)}
    -
    \vec{x}_{k+1\mid k}^{(i)}
    \longrightarrow
    \mathbf P\mathbf H^\top\mathbf\Gamma^{-1}
    \left(
        \mathrm{d}\vec{y}
        -
        \mathbf H\vec{x}^{(i)}\mathrm{d}t
        -
        \mathbf\Gamma^{1/2}\mathrm{d}\vec{W}^{(i)}
    \right).
\end{equation}
Adding Eq.~(\ref{eq:enkf-analysis-limit}) to the model forecast dynamics in Eq.~(\ref{eq:signal_process}) yields the stochastic EnKBF formulation Eq.~(\ref{eq:enkbf-sde-linear}). Applying the derivative-free approximations Eq.~(\ref{eq:abb_tl})--Eq.~(\ref{eq:abb_aj}) then gives Eq.~(\ref{eq:enkbf-sde-stc}).

For the smoother, we start from the discrete-time fixed-interval ensemble RTS smoother. The ensemble-based RTS smoother for Eq.~(\ref{eq:EM-dis_appro}) reads as \citep{evensen2000ensemble}
\begin{equation}\label{eq:EnRTS-update-linear}
    \vec{x}_{k\mid K}^{(i)} =
    \vec{x}_{k\mid k}^{(i)} +
    \mathbf C_k\Big(\vec{x}_{k+1\mid K}^{(i)} -\vec{x}_{k+1\mid k}^{(i)}\Big),
    \qquad i=1,\dots,m,
\end{equation}
where the smoothing gain is given by
\begin{equation}\label{eq:EnRTS-gain}
    \mathbf C_k := \mathbf P_{k|k} (\mathbf I+\tau \mathbf F_k^\top)\big(\mathbf P_{k+1\mid k}\big)^{-1},
\end{equation}
and the Kalman filter prior covariance by \citep{asch2016data,kalman1960new}
\begin{equation*}
    \mathbf P_{k+1\mid k}=\mathbf P_{k\mid k}+\tau\left(\mathbf F_{k} \mathbf P_{k\mid k}+\mathbf P_{k\mid k} \mathbf F_{k}^\top+\mathbf \Sigma_{k}\right) +O(\tau^2).
\end{equation*}
Hence, its inverse admits the following Neumann series expansion for sufficiently small $\tau$ \citep{horn1985matrix,folland1999real}:
    \begin{equation}\label{eq:Neumann}
        \big(\mathbf P_{k+1\mid k})^{-1}
        =( \mathbf P_{k\mid k})^{-1}-\tau ( \mathbf P_{k\mid k})^{-1}\left(\mathbf F_k\mathbf P_{k\mid k}+\mathbf P_{k\mid k}\mathbf F_k^\top+\mathbf \Sigma_{k}\right) (\mathbf P_{k\mid k})^{-1}+ O(\tau^2).
    \end{equation}
Plugging Eq.~(\ref{eq:Neumann}) into Eq.~(\ref{eq:EnRTS-gain}) and expanding again, we obtain
    \begin{equation}\label{eq:C-1}
        (\mathbf C_k)^{-1}=\mathbf I+\tau\big(\mathbf F_k+\mathbf \Sigma_{k}  ( \mathbf P_{k\mid k})^{-1}\big)+ O(\tau^2).
    \end{equation}
Denote the smoother correction for the $i$-th ensemble member by $\vec{e}_{k}^{(i)} := \vec{x}_{k|K}^{(i)} - \vec{x}_{k|k}^{(i)}$. Rearranging Eq.~(\ref{eq:EnRTS-update-linear}) gives
    \begin{equation}\label{eq:e}
        \vec{e}_{k+1}^{(i)}=(\mathbf C_k)^{-1}\vec{e}_{k}^{(i)} -  \left(\vec{x}^{(i)}_{k+1\mid k+1}-\vec{x}^{(i)}_{k+1\mid k}\right).
    \end{equation}
At the formal limit \(\tau\rightarrow0\), we rewrite Eq.~(\ref{eq:enkf-analysis-limit}) as
\[
    \vec{x}^{(i)}_{k+1\mid k+1}-\vec{x}^{(i)}_{k+1\mid k}\rightarrow \mathbf P\mathbf H^\top\mathbf\Gamma^{-1}\mathrm{d} \vec{I}^{(i)},
\]
where \(\mathrm{d} \vec{I}^{(i)} := \mathrm{d} \vec{y} - \mathrm{d} \vec{y}^{(i)}\) denotes the innovation process, with simulated observation $\mathrm{d} \vec{y}^{(i)}$ from Eq.~(\ref{eq:sim-data}).
Omitting now the subscripts \((\cdot)_{k+1}\) or \((\cdot)_{k+1|k+1}\)  in the continuous-time notation,
taking the formal limit $\tau\to 0$ and combining Eq.~(\ref{eq:C-1}) and Eq.~(\ref{eq:e}) yields \citep{law2015data}
    \begin{equation}\label{eq:KBS_err}
        \mathrm{d} \vec{e}^{(i)} = \big(\mathbf F+\mathbf \Sigma \mathbf P^{-1}\big)\vec{e}^{(i)}\mathrm{d} t-\mathbf P\mathbf H^\top\mathbf\Gamma^{-1}\mathrm{d} \vec{I}^{(i)}.
    \end{equation}
Here, \(\vec{e}^{(i)} = \vec{x}_{\mathrm{s}}^{(i)}-\vec{x}_{\mathrm{f}}^{(i)}\).
Adding Eq.~(\ref{eq:KBS_err}) to the filtering member equation
in the linear setting yields a smoother with drift
\(\mathbf F\vec x_{\mathrm s}^{(i)}\).
Replacing this linear drift by
\(\vec f(\vec x_{\mathrm s}^{(i)},t)\) gives the nonlinear
extension Eq.~(\ref{eq:enkbs_member_sde}):
    \begin{equation}\label{eq:enkbs_member_sde_F}
        \mathrm{d} \vec{x}_{\mathrm{s}}^{(i)}
        = \vec{f} \big(\vec{x}_{\mathrm{s}}^{(i)},t\big)\mathrm{d} t+\mathbf \Sigma^{1/2}\mathrm{d} \vec{B}^{(i)}
        + \mathbf\Sigma \mathbf P^{-1}_{\mathrm{f}}\big(\vec{x}_{\mathrm{s}}^{(i)}-\vec{x}_{\mathrm{f}}^{(i)}\big)\mathrm{d} t.
    \end{equation}
Here, $\mathbf P=\mathbf P_{\mathrm{f}}$ is the empirical filtering covariance in the continuous-time notation ($\mathbf P_{k\mid k}\to\mathbf P_{\mathrm{f}}, \tau\to0$). Note that since the smoother dynamics $\mathrm{d} \vec{x}_{\mathrm{s}}^{(i)}$ are obtained by directly adding the correction $\mathrm{d} \vec{e}^{(i)}$ to the forward filter dynamics $\mathrm{d} \vec{x}_{\mathrm{f}}^{(i)}$, the same Wiener process path $\vec{B}^{(i)}$ from the filter is naturally preserved in the smoother dynamics Eq.~(\ref{eq:enkbs_member_sde_F}).

\section{Moment Consistency of the EnKBS in the Linear--Gaussian Regime}
\label{app:moment_consistency}

In this appendix, we demonstrate that the EnKBS equation Eq.~(\ref{eq:enkbs_member_sde}) recovers the classical continuous-time RTS smoothing mean and covariance equations for linear--Gaussian systems \citep{rauch1965maximum, sarkka2023bayesian}.  We note that
all empirical covariance matrices are stochastic objects that depend on the observation path $\vec{y}_{[0,T]}$ and the driving noises $\vec{B}^{(i)}$, $\vec{W}^{(i)}$ for finite ensemble sizes. The analysis below is therefore carried
out in a mean-field limit $m\to\infty$, where empirical moments are identical to the true conditional moments \citep{calvello2025ensemble}. Accordingly, we suppress the ensemble-member index \(i\) throughout
this appendix. 
Throughout this appendix, the observation path $\vec y_{[0,T]}$ is regarded
as fixed, and expectations are taken over the randomness of the
filtering and smoothing processes over $[0,T]$,
including the Wiener noises $\vec B$ and $\vec W$.

Consider the continuous-time linear--Gaussian setting where the drift and observation maps in Eq.~(\ref{eq:general_system}) are linear: $\vec{f} = \mathbf F\vec{x}$ and $\vec{h}= \mathbf H \vec{x}$. 
Let $\vec{\hat x}_{\mathrm{s}} := \vec{x}_{\mathrm{s}} - \overline{\vec{x}}_{\mathrm{s}}$ and $\vec{\hat x}_{\mathrm{f}} := \vec{x}_{\mathrm{f}} - \overline{\vec{x}}_{\mathrm{f}}$ denote the smoother and filter anomalies, respectively. By averaging over ensemble members \(i\) in the EnKBS member equation Eq.~(\ref{eq:enkbs_member_sde}) and taking a mean-field (infinite ensemble) limit \(m\rightarrow\infty\), we recover the exact mean equation of the RTS smoother:
\begin{equation}\label{eq:enkbs_mean}
        \mathrm{d} \overline{\vec{x}}_{\mathrm{s}}
        = \mathbf F\overline{\vec{x}}_{\mathrm{s}}\mathrm{d} t+ \mathbf\Sigma \mathbf P^{-1}_{\mathrm{f}}\big(\overline{\vec{x}}_{\mathrm{s}}-\overline{\vec{x}}_{\mathrm{f}}\big)\mathrm{d} t.
\end{equation} 
Subtracting Eq.~(\ref{eq:enkbs_mean}) from the EnKBS member equation Eq.~(\ref{eq:enkbs_member_sde}) yields the anomaly dynamics:
\begin{equation}\label{eq:app_anomaly_sde}
    \mathrm{d} \vec{\hat x}_{\mathrm{s}}
    = \left(( \mathbf F + \mathbf\Sigma \mathbf P^{-1}_{\mathrm{f}})\vec{\hat x}_{\mathrm{s}}-\mathbf\Sigma \mathbf P^{-1}_{\mathrm{f}}\vec{\hat x}_{\mathrm{f}}\right)\mathrm{d} t + \mathbf \Sigma^{1/2}\mathrm{d} \vec{B}.
\end{equation}

For the stochastic EnKBF, averaging Eq.~(\ref{eq:enkbf-sde-stc}) over the ensemble gives
\begin{equation}\label{eq:app_filter_mean}
    \mathrm{d}\overline{\vec{x}}_{\mathrm f}
    =
    \mathbf F\overline{\vec{x}}_{\mathrm f}\mathrm{d}t
    +
    \mathbf P_{\mathrm f}\mathbf H^\top\mathbf\Gamma^{-1}
    \left(
        \mathrm{d}\vec y
        -
        \mathbf H\overline{\vec{x}}_{\mathrm f}\mathrm{d}t
    \right),
\end{equation}
which is the standard Kalman--Bucy filtering mean equation
\citep{bergemann2012Ensemble,jazwinski2007stochastic}. Its filtering
anomalies analogously satisfy
\begin{equation}\label{eq:app_enkf_anomaly}
    \mathrm{d}\vec{\hat x}_{\mathrm{f}}
    =
    (\mathbf F-\mathbf K\mathbf H)
    \vec{\hat x}_{\mathrm{f}}\mathrm{d}t
    +
    \mathbf\Sigma^{1/2}\mathrm{d}\vec B
    -
    \mathbf K\mathbf\Gamma^{1/2}\mathrm{d}\vec W,
\end{equation}
where
\(
    \mathbf K
    =
    \mathbf P_{\mathrm f}\mathbf H^\top\mathbf\Gamma^{-1}
\)
is the Kalman--Bucy filter gain.
It is useful to write the filtering anomaly dynamics more generally as
\begin{equation}\label{eq:app_filter_general}
    \mathrm{d}\vec{\hat x}_{\mathrm f}
    =
    \mathbf A_{\mathrm f}
    \vec{\hat x}_{\mathrm f}\mathrm{d}t
    +
    \mathbf\Sigma^{1/2}\mathrm{d}\vec B
    +
    \mathbf L\mathrm{d}\vec W.
\end{equation}
For the stochastic EnKBF,
\begin{equation}\label{eq:app_enkbf_params}
        \mathbf A_{\mathrm f}
    =
    \mathbf F-\mathbf K\mathbf H,
    \qquad
    \mathbf L
    =
    -\mathbf K\mathbf\Gamma^{1/2}.
\end{equation}
Applying It\^o's product rule \(\mathrm d(XX^\top)
=(\mathrm dX)X^\top
+X(\mathrm dX)^\top
+(\mathrm dX)(\mathrm dX)^\top\) to
$\vec{\hat x}_{\mathrm f}
(\vec{\hat x}_{\mathrm f})^\top$
and Eq.~(\ref{eq:app_filter_general}) yields
\citep{oksendal2003Stochastic,sarkka2019applied}:
\begin{equation}\label{eq:app_filter_cov_general}
    \dot{\mathbf P}_{\mathrm f}
    =
    \mathbf A_{\mathrm f}\mathbf P_{\mathrm f}
    +
    \mathbf P_{\mathrm f}\mathbf A_{\mathrm f}^{\top}
    +
    \mathbf\Sigma
    +
    \mathbf L\mathbf L^\top.
\end{equation}
Substituting Eq.~(\ref{eq:app_enkbf_params}) into
Eq.~(\ref{eq:app_filter_cov_general}) gives
\begin{align}
\dot{\mathbf P}_{\mathrm f}
={}&
(\mathbf F-\mathbf K\mathbf H)\mathbf P_{\mathrm f}
+
\mathbf P_{\mathrm f}
(\mathbf F-\mathbf K\mathbf H)^\top
+
\mathbf\Sigma
+
\mathbf K\mathbf\Gamma\mathbf K^\top
\notag\\
={}&
\mathbf F\mathbf P_{\mathrm f}
+
\mathbf P_{\mathrm f}\mathbf F^\top
+
\mathbf\Sigma
-
\mathbf P_{\mathrm f}\mathbf H^\top
\mathbf\Gamma^{-1}\mathbf H\mathbf P_{\mathrm f},
\label{eq:app_kb_riccati}
\end{align}
which is the standard Kalman--Bucy filtering covariance equation \citep{kalman1961new, jazwinski2007stochastic}.

We note that the forward--backward construction considered here
differs from conventional forward--backward SDEs
\citep{ma1994solving}. Unlike the filtering anomaly, the fixed-interval
smoothing anomaly
\(\vec{\hat x}_{\mathrm s}(t)\) is obtained by propagating
information backward from the terminal time \(T\). At
time \(t\), it depends on Wiener increments over both \([0,t]\) and
\((t,T]\), rather than only on the noise history up to \(t\). Thus, \(\vec{\hat x}_{\mathrm s}(t)\) does not satisfy the forward It\^o-process form as required by the standard It\^o product rule \citep{oksendal2003Stochastic,sarkka2019applied}.
Consequently, the standard forward It\^o product-rule calculation used
to obtain Eq.~(\ref{eq:app_filter_cov_general}) is not directly
applicable to
\(\vec{\hat x}_{\mathrm s}
(\vec{\hat x}_{\mathrm s})^\top\).
We instead derive the smoothing covariance below from its integral
representation.

We now return to the smoother dynamics. Define
\begin{equation}\label{eq:app_As_e}
    \mathbf A_{\mathrm s}
    :=
    \mathbf F+\mathbf\Sigma\mathbf P_{\mathrm f}^{-1},
    \qquad
    \mathbf D
    :=
    \mathbf F-\mathbf A_{\mathrm f},
    \qquad
    \vec e
    :=
    \vec{\hat x}_{\mathrm s}
    -
    \vec{\hat x}_{\mathrm f}.
\end{equation}
Here, $\mathbf A_{\mathrm f}$ and $\mathbf A_{\mathrm s}$ are the linear
generators of the filter-anomaly and smoother-anomaly dynamics, respectively.
Subtracting Eq.~(\ref{eq:app_filter_general}) from
Eq.~(\ref{eq:app_anomaly_sde}), and recalling that the filter and smoother
use the same realization of the system noise
$\mathrm{d}\vec B$, gives
\begin{equation}\label{eq:app_error_terminal}
    \mathrm{d}\vec e
    =
    \mathbf A_{\mathrm s}\vec e\mathrm{d}t
    +
    \mathbf D\vec{\hat x}_{\mathrm f}\mathrm{d}t
    -
    \mathbf L\mathrm{d}\vec W,
    \qquad
    \vec e(T)=\vec 0.
\end{equation}
The shared system-noise increment $\mathrm{d}\vec B$ cancels exactly in this difference.

Let $\mathbf\Phi_{\mathrm f}(t,u)$ and
$\mathbf\Phi_{\mathrm s}(t,u)$ denote the fundamental matrices
generated by $\mathbf A_{\mathrm f}$ and $\mathbf A_{\mathrm s}$,
respectively, such that
\begin{equation}\label{eq:app_transition}
    \frac{\partial}{\partial t}\mathbf\Phi_j(t,u)
    =
    \mathbf A_j(t)\mathbf\Phi_j(t,u),
    \qquad
    \mathbf\Phi_j(u,u)=\mathbf I,
    \qquad
    j\in\{\mathrm f,\mathrm s\}.
\end{equation}
These satisfy the composition or semigroup property
\begin{equation}\label{eq:app_transition_semigroup}
    \mathbf\Phi_j(t,s)\mathbf\Phi_j(s,u)
    =
    \mathbf\Phi_j(t,u).
\end{equation}
Solving the terminal-value equation
Eq.~(\ref{eq:app_error_terminal}) by variation of constants gives \citep{chicone2024ordinary}
\begin{equation}\label{eq:app_error_solution}
    \vec e(t)
    =
    -\int_t^T
    \mathbf\Phi_{\mathrm s}(t,u)\mathbf D(u)
    \vec{\hat x}_{\mathrm f}(u)\mathrm{d}u
    +
    \int_t^T
    \mathbf\Phi_{\mathrm s}(t,u)\mathbf L(u)
    \mathrm{d}\vec W_u.
\end{equation}
Likewise, for $u\geq t$, solving the forward filtering anomaly Eq.~(\ref{eq:app_filter_general}) gives
\begin{equation}\label{eq:app_filter_solution}
\begin{aligned}
    \vec{\hat x}_{\mathrm f}(u)
    ={}
    \mathbf\Phi_{\mathrm f}(u,t)
    \vec{\hat x}_{\mathrm f}(t)
    +
    \int_t^u
    \mathbf\Phi_{\mathrm f}(u,r)
    \mathbf\Sigma^{1/2}\mathrm{d}\vec B_r
    +
    \int_t^u
    \mathbf\Phi_{\mathrm f}(u,r)
    \mathbf L(r)\mathrm{d}\vec W_r.
\end{aligned}
\end{equation}

For notational convenience, introduce the auxiliary matrices
\begin{equation}\label{eq:app_CM}
    \mathbf C(t)
    :=
    \int_t^T
    \mathbf\Phi_{\mathrm s}(t,u)
    \mathbf D(u)
    \mathbf\Phi_{\mathrm f}(u,t)\mathrm{d}u,
    \qquad
    \mathbf M(t):=\mathbf I-\mathbf C(t).
\end{equation}
The matrix $\mathbf C(t)$ accumulates the effect of the filter-update
coupling $\mathbf D$ over the future interval $[t,T]$, with the filter
dynamics propagated forward and the smoother dynamics propagated backward.
Substituting Eq.~(\ref{eq:app_filter_solution}) into
Eq.~(\ref{eq:app_error_solution}), formally interchanging the order of integration,
\begin{align}
    \int_t^T
    \mathbf\Phi_{\mathrm s}(t,u)\mathbf D(u)
    \left[
    \int_t^u
    \mathbf\Phi_{\mathrm f}(u,r)
    \mathbf\Sigma^{1/2}
    \mathrm{d}\vec B_r
    \right]
    \mathrm{d}u
    &=
    \int_t^T
    \left[
    \int_r^T
    \mathbf\Phi_{\mathrm s}(t,u)
    \mathbf D(u)
    \mathbf\Phi_{\mathrm f}(u,r)
    \mathrm{d}u
    \right]
    \mathbf\Sigma^{1/2}
    \mathrm{d}\vec B_r
    \notag\\
    &=
    \int_t^T
    \mathbf\Phi_{\mathrm s}(t,r)
    \mathbf C(r)
    \mathbf\Sigma^{1/2}
    \mathrm{d}\vec B_r,
    \label{eq:app_fubini_B}
\end{align}
and similarly,
\begin{equation}\label{eq:app_fubini_W}
    \int_t^T
    \mathbf\Phi_{\mathrm s}(t,u)\mathbf D(u)
    \left[
    \int_t^u
    \mathbf\Phi_{\mathrm f}(u,r)
    \mathbf L(r)
    \mathrm{d}\vec W_r
    \right]
    \mathrm{d}u
    =
    \int_t^T
    \mathbf\Phi_{\mathrm s}(t,r)
    \mathbf C(r)\mathbf L(r)
    \mathrm{d}\vec W_r,
\end{equation}
we obtain
\begin{equation}\label{eq:app_error_representation}
\begin{aligned}
    \vec e(t)
    ={}
    -\mathbf C(t)\vec{\hat x}_{\mathrm f}(t)
    -
    \int_t^T
    \mathbf\Phi_{\mathrm s}(t,r)
    \mathbf C(r)\mathbf\Sigma^{1/2}
    \mathrm{d}\vec B_r
    +
    \int_t^T
    \mathbf\Phi_{\mathrm s}(t,r)
    \mathbf M(r)\mathbf L(r)
    \mathrm{d}\vec W_r.
\end{aligned}
\end{equation}
Since
\(\vec{\hat x}_{\mathrm s}
    =
    \vec{\hat x}_{\mathrm f}
    +
    \vec e,\)
and $\mathbf M=\mathbf I-\mathbf C$,
Eq.~(\ref{eq:app_error_representation}) gives
\begin{equation}\label{eq:app_smoother_representation}
\begin{aligned}
    \vec{\hat x}_{\mathrm s}(t)
    ={}
    \mathbf M(t)
    \vec{\hat x}_{\mathrm f}(t)
    -
    \int_t^T
    \mathbf\Phi_{\mathrm s}(t,r)
    \mathbf C(r)\mathbf\Sigma^{1/2}
    \mathrm{d}\vec B_r+
    \int_t^T
    \mathbf\Phi_{\mathrm s}(t,r)
    \mathbf M(r)\mathbf L(r)
    \mathrm{d}\vec W_r.
\end{aligned}
\end{equation}

We now define the smoother posterior covariance by
\[
    \mathbf P_{\mathrm s}(t)
    :=
    \mathbb E\left[
    \vec{\hat x}_{\mathrm s}(t)
    \big(\vec{\hat x}_{\mathrm s}(t)\big)^\top
    \right].
\]
In the linear--Gaussian mean-field limit, $\mathbf P_{\mathrm f}(t)$
converges to the conditional filtering covariance and is therefore
deterministic for a fixed observation trajectory. Consequently,
$\mathbf A_{\mathrm f}$, $\mathbf A_{\mathrm s}$, $\mathbf D$, $\mathbf L$,
$\mathbf\Phi_{\mathrm f}$, $\mathbf\Phi_{\mathrm s}$, $\mathbf C$, and
$\mathbf M$ are deterministic functions of time.
It follows that the first term in
Eq.~(\ref{eq:app_smoother_representation}) depends on the filtering
noises up to time $t$, whereas the two stochastic-integral terms
depend only on Wiener increments over $(t,T]$.
Therefore, these
two integrals are uncorrelated with
$\vec{\hat x}_{\mathrm f}(t)$. Using the independence
of $\vec B$ and $\vec W$, and \(\mathrm{d}\vec B_r \left(\mathrm{d} \vec B_r\right)^\top = \mathbf I  \mathrm{d}r\), \(\mathrm{d}\vec W_r \left(\mathrm{d}\vec W_r\right)^\top = \mathbf I  \mathrm{d}r\),
Eq.~(\ref{eq:app_smoother_representation}) gives
\begin{equation}\label{eq:app_Ps_representation}
    \mathbf P_{\mathrm s}
    =
    \mathbf M\mathbf P_{\mathrm f}\mathbf M^\top
    +
    \mathbf N,
\end{equation}
where
\begin{equation}\label{eq:app_N}
\begin{aligned}
    \mathbf N(t)
    :=
    \int_t^T
    \mathbf\Phi_{\mathrm s}(t,r)
    \Big[
        \mathbf C(r)\mathbf\Sigma\mathbf C(r)^\top
        +
        \mathbf M(r)\mathbf L(r)\mathbf L(r)^\top
        \mathbf M(r)^\top
    \Big]
    \mathbf\Phi_{\mathrm s}(t,r)^\top
    \mathrm{d}r.
\end{aligned}
\end{equation}

It remains to differentiate Eq.~(\ref{eq:app_Ps_representation}).
From Eq.~(\ref{eq:app_CM}), using Eq.~(\ref{eq:app_transition}), we obtain
\begin{equation}\label{eq:app_C_ode}
    \dot{\mathbf C}
    =
    -\mathbf D
    +
    \mathbf A_{\mathrm s}\mathbf C
    -
    \mathbf C\mathbf A_{\mathrm f}.
\end{equation}
Since
$\mathbf M=\mathbf I-\mathbf C$, and
$\mathbf A_{\mathrm s}
=\mathbf F+\mathbf\Sigma\mathbf P_{\mathrm f}^{-1}$,
Eq.~(\ref{eq:app_C_ode}) is equivalent to
\begin{equation}\label{eq:app_M_ode}
    \dot{\mathbf M}
    =
    \mathbf A_{\mathrm s}\mathbf M
    -
    \mathbf M\mathbf A_{\mathrm f}
    -
    \mathbf\Sigma\mathbf P_{\mathrm f}^{-1},
    \qquad
    \mathbf M(T)=\mathbf I.
\end{equation}
Differentiating Eq.~(\ref{eq:app_N}) with respect to its lower limit $t$
similarly gives
\begin{equation}\label{eq:app_N_ode}
\begin{aligned}
    \dot{\mathbf N}
    =
    \mathbf A_{\mathrm s}\mathbf N
    +
    \mathbf N\mathbf A_{\mathrm s}^\top
    -
    \mathbf C\mathbf\Sigma\mathbf C^\top
    -
    \mathbf M\mathbf L\mathbf L^\top\mathbf M^\top,
    \qquad
    \mathbf N(T)=\mathbf 0.
\end{aligned}
\end{equation}
Differentiating Eq.~(\ref{eq:app_Ps_representation}) yields
\begin{equation}\label{eq:app_Ps_product}
    \dot{\mathbf P}_{\mathrm s}
    =
    \dot{\mathbf M}\mathbf P_{\mathrm f}\mathbf M^\top
    +
    \mathbf M\dot{\mathbf P}_{\mathrm f}\mathbf M^\top
    +
    \mathbf M\mathbf P_{\mathrm f}\dot{\mathbf M}^\top
    +
    \dot{\mathbf N}.
\end{equation}
Substituting Eqs.~(\ref{eq:app_filter_cov_general}),
(\ref{eq:app_M_ode}), and (\ref{eq:app_N_ode}) into
Eq.~(\ref{eq:app_Ps_product}), all terms involving
$\mathbf A_{\mathrm f}$ and $\mathbf L\mathbf L^\top$ cancel. Using
Eq.~(\ref{eq:app_Ps_representation}), the remaining terms reduce to
\begin{equation}\label{eq:app_before_final_cancel}
    \dot{\mathbf P}_{\mathrm s}
    =
    \mathbf A_{\mathrm s}\mathbf P_{\mathrm s}
    +
    \mathbf P_{\mathrm s}\mathbf A_{\mathrm s}^\top
    -
    \mathbf\Sigma\mathbf M^\top
    -
    \mathbf M\mathbf\Sigma
    +
    \mathbf M\mathbf\Sigma\mathbf M^\top
    -
    \mathbf C\mathbf\Sigma\mathbf C^\top.
\end{equation}
Finally, since $\mathbf C=\mathbf I-\mathbf M$,
\begin{equation}
    -\mathbf C\mathbf\Sigma\mathbf C^\top
    =
    -\mathbf\Sigma
    +
    \mathbf M\mathbf\Sigma
    +
    \mathbf\Sigma\mathbf M^\top
    -
    \mathbf M\mathbf\Sigma\mathbf M^\top.
\end{equation}
All $\mathbf M$-dependent terms in
Eq.~(\ref{eq:app_before_final_cancel}) also cancel, giving
\begin{equation}\label{eq:rts_cov}
    \dot{\mathbf P}_{\mathrm s}
    =
    \left(
    \mathbf F+\mathbf\Sigma\mathbf P_{\mathrm f}^{-1}
    \right)\mathbf P_{\mathrm s}
    +
    \mathbf P_{\mathrm s}
    \left(
    \mathbf F+\mathbf\Sigma\mathbf P_{\mathrm f}^{-1}
    \right)^\top
    -
    \mathbf\Sigma,
    \qquad
    \mathbf P_{\mathrm s}(T)=\mathbf P_{\mathrm f}(T).
\end{equation}
This is the exact continuous-time RTS smoother covariance equation
\citep{rauch1965maximum,sarkka2023bayesian}. Together with
Eq.~(\ref{eq:enkbs_mean}), this shows that the EnKBS recovers both the
RTS smoothing mean and covariance equations in the linear--Gaussian
mean-field limit.

Importantly, this conclusion is not specific to the stochastic EnKBF
realization. To see this, consider the continuous-time deterministic
EnKBF
\citep{sakov2008deterministic,bergemann2012Ensemble}:
\begin{equation}\label{eq:enkbf-sde-deter}
    \mathrm{d}\vec{x}_{\mathrm f}
    =
    \mathbf F\vec{x}_{\mathrm f}\mathrm{d}t
    +
    \mathbf\Sigma^{1/2}\mathrm{d}\vec B
    +
    \mathbf P_{\mathrm f}\mathbf H^\top\mathbf\Gamma^{-1}
    \left(
        \mathrm{d}\vec y
        -
        \frac{
        \mathbf H\vec{x}_{\mathrm f}
        +
        \mathbf H\overline{\vec{x}}_{\mathrm f}
        }{2}\mathrm{d}t
    \right).
\end{equation}
Averaging Eq.~(\ref{eq:enkbf-sde-deter}) gives the same Kalman--Bucy
filtering mean equation as the stochastic EnKBF. Its anomalies satisfy
\begin{equation}\label{eq:denkbf_anomaly}
    \mathrm{d}\vec{\hat x}_{\mathrm f}
    =
    \left(
        \mathbf F-\frac12\mathbf K\mathbf H
    \right)
    \vec{\hat x}_{\mathrm f}\mathrm{d}t
    +
    \mathbf\Sigma^{1/2}\mathrm{d}\vec B.
\end{equation}
Thus, in the general representation
Eq.~(\ref{eq:app_filter_general}),
\[
    \mathbf A_{\mathrm f}
    =
    \mathbf F-\frac12\mathbf K\mathbf H,
    \qquad
    \mathbf L=\mathbf 0.
\]
In the linear--Gaussian mean-field limit, the arguments in this appendix
apply in the same manner to this deterministic forward formulation. Note that the cancellations leading from
Eq.~(\ref{eq:app_Ps_product}) to Eq.~(\ref{eq:rts_cov})
are independent of the specific choices of
\(\mathbf A_{\mathrm f}\) and \(\mathbf L\) in the forward EnKBF
formulation and therefore remain valid for the deterministic EnKBF.
For these choices of \(\mathbf A_{\mathrm f}\) and \(\mathbf L\), the
deterministic EnKBF recovers the same Kalman--Bucy filtering mean and
covariance equations Eqs.~(\ref{eq:app_filter_mean}) and
(\ref{eq:app_kb_riccati}). In addition, when either the stochastic or
deterministic EnKBF formulation is used in the forward pass, the
corresponding EnKBS recovers the continuous-time RTS smoothing mean and
covariance equations Eqs.~(\ref{eq:enkbs_mean}) and
(\ref{eq:rts_cov}).

We have also tested the deterministic EnKBF
Eq.~(\ref{eq:enkbf-sde-deter}) as the forward pass in our framework. It
yields results similar to those of the stochastic EnKBF
Eq.~(\ref{eq:enkbf-sde-stc}) across all numerical experiments presented
in this paper (not shown). In addition, for either choice of the forward
EnKBF, the smoother consistently improves upon the filter by reducing
both the estimation error and the ensemble spread (with behavior similar
to that shown in Figure~\ref{fig:good_x_to_y_err_m50}).

\noappendix       







\authorcontribution{NC and MA conceived the project, with NC designing the numerical experiments. ZJ developed the code, performed experiments, and prepared the manuscript draft and figures, in regular discussion with NC and MA. SR reviewed an earlier version and provided critical feedback that motivated the analysis in Appendix~\ref{app:moment_consistency}. NC supervised the project, and all authors contributed to manuscript revision.} 

\competinginterests{NC is a member of the editorial board of \textit{Nonlinear Processes in Geophysics}. The authors have no other competing interests to declare.} 


\begin{acknowledgements}
    The research of NC is funded by the Army Research Office W911NF-23-1-0118. ZJ and MA are partially supported as research assistants under this grant. MA is also supported in part by NSF Award DMS-2023239. The work of SR has been partially funded by Deutsche Forschungsgemeinschaft (DFG) -- Project-ID 318763901 -- SFB1294.
\end{acknowledgements}







\bibliographystyle{copernicus}
\bibliography{refs.bib}

\end{document}